\documentclass[12pt]{article}
\usepackage[cp1251]{inputenc}
\usepackage[russian, english]{babel}
\usepackage{amsfonts,amsmath}
\usepackage{graphicx} 
\usepackage[labelsep=period]{caption}
\usepackage{tikz-cd}
\usepackage{xcolor,float}

\usepackage[left=1.5cm, right=1.5cm, top=2cm, bottom=2cm, headsep=0.7cm, footskip=1cm
]{geometry}
\newtheorem{theorem}{Theorem}
\newtheorem{lemma}{Lemma}

\newtheorem{remark}{Remark}

\begin{document}

\title{Billiard trajectories inside Cones}

\author {Andrey E. Mironov and Siyao Yin}

\date{}
\maketitle

\sloppy
\begin{abstract}
  In \cite{MY}, it was proved that every billiard trajectory inside a $C^3$ convex cone has a finite number of reflections. 
  Here, by a $C^3$ convex cone, we mean a cone whose section with some hyperplane is a strictly convex closed $C^3$ submanifold of the hyperplane with nondegenerate second fundamental form.
    In this paper, 
  we prove the existence of $C^2$ convex cones admitting billiard trajectories with infinitely many reflections in finite time.
  We also estimate the number of reflections of billiard trajectories in elliptic cones in $\mathbb{R}^3$
  using two first integrals.
\end{abstract}

\tableofcontents
  
\section{Introduction and Main Results}
\begin{figure}[h]
    \begin{center}
    \includegraphics[scale=0.47]{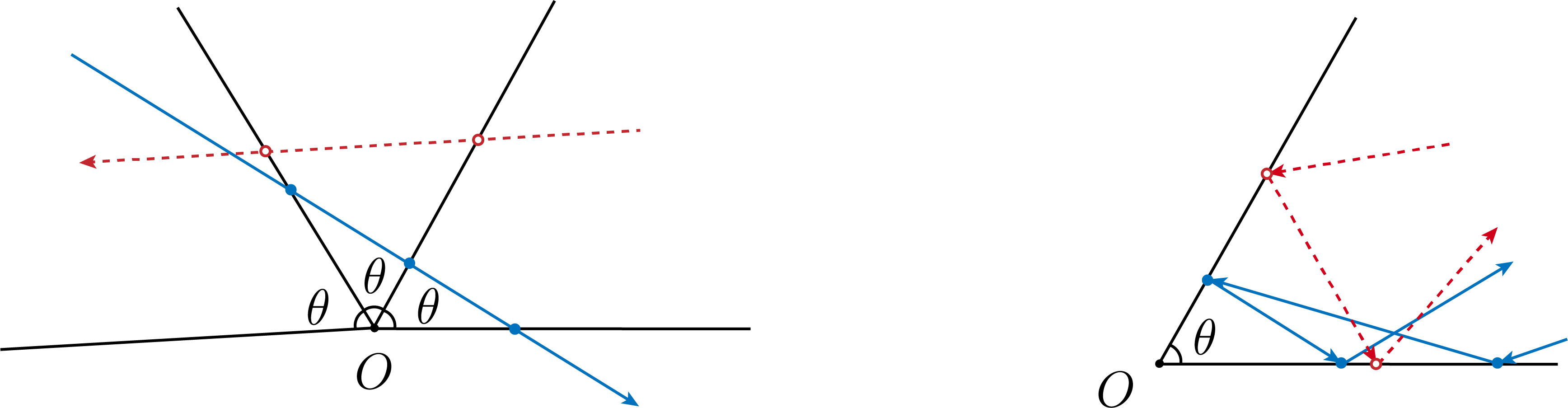}
    \end{center}
    \caption{For an angle $\theta$ in $\mathbb{R}^2$, a trajectory has $n$ or $n-1$ reflections, where $\frac{\pi}{n} \leq \theta < \frac{\pi}{n-1}$.
      }
    \label{fig:angle-theta}
  \end{figure}

The Birkhoff billiard is an important dynamical system that has inspired many remarkable conjectures and results (see e.g., \cite{Bir}--\cite{Bol}).
In this paper we continue to study Birkhoff billiards inside cones $K \subset \mathbb{R}^n$ started in \cite{MY}. 
In the simplest case $n=2$, the cone $K$ reduces to an angle $\theta$, inside which 
the billiard dynamics is well understood.
By repeatedly reflecting the angle $\theta$ along its sides, one can cover the plane with copies of the angle. The billiard trajectory within $\theta$ then becomes a straight line in the plane with each intersection corresponding to a reflection (see Fig. \ref{fig:angle-theta}).  
Therefore the number of reflections is at most $\left\lceil \frac{\pi}{\theta} \right\rceil$, where $\left\lceil x \right\rceil$ is the smallest integer greater than or equal to $x$ (see e.g., \cite{Tab}, \cite{Gal}).  

Sinai \cite{Sin} studied billiard dynamics inside a polyhedral angle and proved that, as in the two-dimensional case, the number of reflections of a billiard trajectory is finite.  
Moreover, a uniform bound on the number of reflections for all trajectories inside a fixed polyhedral angle was established (see also \cite{Sev}).

In \cite{MY}, an analog of Sinai's Theorem in the smooth case was proved. More precisely, consider the cone $K$ defined as  
\begin{equation}\label{eq:def-k}
  K = \{ tp \mid p \in \gamma, t > 0 \} \subset \mathbb{R}^n,
\end{equation}
where $\gamma \subset \mathcal{P}:= \{ x \in \mathbb{R}^n \mid x^n = 1 \}$ is a $C^3$-smooth, strictly convex, closed submanifold of $\mathcal{P}$ with nondegenerate second fundamental form at every point. 
Then any billiard trajectory in $K$ has a finite number of reflections. 
Building on this result, the integrability of the billiard system inside $K$ was established in \cite{MY}, i.e., a set of first integrals was constructed on the phase space of the system, such that each value of this set of first integrals uniquely determines a trajectory inside $K$.

We say that $K$, defined by (\ref{eq:def-k}), is $C^2$ smooth if $\gamma$ is $C^2$ smooth. The main result of this paper is the following.
\begin{theorem}\label{thm:c2}
  There exist $C^2$-smooth  convex cones with billiard trajectories having infinitely many reflections in finite time.
\end{theorem}

For billiards in ${\mathbb R}^2$, Halpern \cite{Hal} showed that $C^2$-smooth, strictly convex tables can admit trajectories with infinitely many reflections in finite time, but such trajectories cannot exist when the boundary is $C^3$-smooth. Gruber \cite{Gru2} later extended this result to ${\mathbb R}^n$, proving the non-existence of such trajectories for $C^3$-smooth, strictly convex hypersurfaces. Theorem \ref{thm:c2} is an analog of Halpern's theorem for closed curves, to the case of cones over  closed hypersurfaces.

It is evident that there is no universal constant bounding the number of reflections in a smooth cone, i.e., an analog of Sinai's uniform estimate does not hold for smooth cones.
This raises the question: under what conditions can the number of reflections of a billiard trajectory be estimated?
For the elliptic cone in ${\mathbb R}^3$, we can estimate the number of reflections when the values of the first integrals are fixed. 
Let $K_e$ be an elliptic cone in $\mathbb{R}^3$ defined by the equation  
\begin{equation}\label{eq:elliptic-cone}
(x^3)^2 = \frac{(x^1)^2}{a^2} + \frac{(x^2)^2}{b^2}, \quad x^3 > 0, \quad a>b>0.
\end{equation}  
Then the Birkhoff billiard inside $K_e$ admits the following first integrals (see \cite{MY} and Lemma \ref{lem:int-ell} below):  
$$
I_1 = 
m_{1,2}^2 + m_{1,3}^2 + m_{2,3}^2,
\quad 
I_2 = 
a^2 m_{2,3}^2 + b^2 m_{1,3}^2 - m_{1,2}^2,
$$
where $m_{i,j} := x^i v^j - x^j v^i$ for $i < j$, $i, j = 1, 2, 3$, and ${v} = (v^1, v^2, v^3)$ is the velocity vector.
One can check that the Poisson bracket of $I_1$ and $I_2$ is $0$,
$$
 \{I_1,I_2\}=\sum_{k=1}^3
 \left(
  \frac{\partial I_1}{\partial x^k}\frac{\partial I_2}{\partial v^k}-
 \frac{\partial I_1}{\partial v^k}\frac{\partial I_2}{\partial x^k}\right)=0.
$$

Our next result is the following.
  
\begin{theorem}\label{thm:elliptic} 
    The number of reflections of a trajectory inside $K_e$ with fixed values of the first integrals
    \begin{equation}\label{eq:i-c}
      I_1 = c_1>0 , \quad I_2 = c_2>0
    \end{equation}
    is bounded by $N_{c_1,c_2}$, where
    $$
    N_{c_1,c_2} = \left\lceil 
    \frac{\pi}
    {\arcsin 
    \frac{2 a b \sqrt{c_1 c_2}}
    {a^2(b^2+1)c_1 + (b^2+1)c_2}
    } 
    \right\rceil.
    $$
  \end{theorem}

  The paper is organized as follows. 
  We prove Theorem \ref{thm:c2} in Section 2, and Theorem \ref{thm:elliptic} in Section 3.

  The authors are grateful to Misha Bialy for valuable discussions and suggestions.

\section{Billiard Trajectories in $C^2$ Cones}
In this section we prove Theorem \ref{thm:c2}.

\subsection{Geometric properties of billiard trajectories}  
We begin by recalling several geometric properties of billiard trajectories inside cones, which are derived in \cite{MY} and will be used in this paper. These properties apply to the cone $K$ over any hypersurface $\gamma \subset \mathcal{P}$, without assuming specific smoothness or convexity conditions on $\gamma$. 

\begin{lemma}[Theorem 1, Lemma 2 in \cite{MY}]\label{lem:int-dis}  
  Let $l \subset \mathbb{R}^n$ be an oriented line with direction vector ${v} = (v^1, \ldots, v^n)$, where $\|{v}\| = 1$. Define the function  
  $$
  I = \sum_{1 \leq i < j \leq n} m_{i,j}^2, 
  $$  
  where $m_{i,j} := x^i v^j - x^j v^i$, $i < j$, $i,j = 1, \ldots, n$, and $x=(x^1, \ldots, x^n) \in l$.  
  Then $I$ is a first integral for the billiard trajectory inside $K$ and represents the square of the distance from $l$ to the origin.
  \end{lemma} 
The distance from $l$ to the origin can equivalently be expressed as (see Equation (2) in \cite{MY}):
\begin{equation}\label{eqdistance}  
\text{dist}(l, O) = \left(\|x\|^2 - \langle x, {v} \rangle^2 \right)^{1/2}.
\end{equation}  

Now, consider a billiard trajectory inside $K$ with infinitely many reflections.  
Let $l_k$ be the sequence of oriented lines associated with this trajectory, where $l_k$ has direction vector ${v}_{k}$. Each $l_k$ intersects $K$ at points $p_k$ and $p_{k+1}$, satisfying  
$
\frac{p_{k+1} - p_k}{\|p_{k+1} - p_k\|} = {v}_{k}.
$
Let $\alpha_k \in (0, \pi)$ denote the angle between ${v}_{k}$ and $p_k$, and $\theta_k \in (0, \pi)$ denote the angle $\angle p_k O p_{k+1}$. Then the following relation holds (see Lemma 3 in \cite{MY}): 
\begin{equation}\label{eq:alpha-theta}  
\alpha_{k+1} = \alpha_k - \theta_k.  
\end{equation}  
In addition, there is a relationship between $\|p_k\|$ and $\sin \alpha_k$, given by (see Equation (8) in \cite{MY}):
\begin{equation}\label{eq:sin-p1}  
\|p_k\| = \frac{R}{\sin \alpha_k},  
\end{equation}  
where $R$ is the distance from the line $l_k$ to the origin. By Lemma \ref{lem:int-dis}, $R$ is constant for all $k$.  
  
\subsection{$C^2$ cones with infinite reflections}
  In this section, we construct a $C^2$ convex cone $K$ that admits trajectories with infinitely many reflections.
These trajectories exhibit two types of limit behavior depending on a parameter 
$a \in (-\frac{\pi}{2}, \frac{\pi}{2}]$: for $-\frac{\pi}{2} < a < \frac{\pi}{2}$, the reflection points converge to a finite limit point; for $a = \frac{\pi}{2}$, the reflection points tend to infinity and the trajectory is asymptotic to a half-line in $K$.

We begin by presenting the construction in $\mathbb{R}^3$.
At the end of this section, we explain how this construction can be generalized to $\mathbb{R}^n$.

To do this, we first construct a sequence of points $p_k$ that will be the vertices of the billiard trajectory inside $K$.
Let $K_0 \subset \mathbb{R}^3$ be the circular cone defined by
$x^3 = \sqrt{(x^1)^2 + (x^2)^2}$, 
and $\gamma_0$ be the 
the intersection of $K_0$ with the plane $\mathcal{P} = \{x \in \mathbb{R}^3 \mid x^3 = 1\}$, 
i.e., $\gamma_0$ is the unit circle $\{(x^1,x^2,1)\mid (x^1)^2 + (x^2)^2 = 1\}$.

Consider a sequence of points
\begin{equation}\label{eq:def-of-pk}
  p_k := (t_k \cos \xi_k, t_k \sin \xi_k, t_k),
\quad \xi_k = k^{-1/2},
\quad t_k >0,
\quad k\geq 1,
\end{equation}
where $t_k$ will be defined later in (\ref{eq:tk_a}).
Given the values of $\xi_k, k\geq 1$, the angle $\theta_k=\angle p_k O p_{k+1}$ is uniquely determined and 
does not depend on $t_k$ and $t_{k+1}$.
Indeed, 
\begin{equation}\label{eq:xi-theta-0}
\cos \theta_k = \left\langle \frac{p_k}{\|p_k\|}, \frac{p_{k+1}}{\|p_{k+1}\|} \right\rangle 
= \frac{\cos \xi_{k}\cos \xi_{k+1} + \sin \xi_k \sin \xi_{k+1} + 1}{2}
= \frac{\cos (\xi_k - \xi_{k+1}) + 1 }{2} 
.
\end{equation}
Since $\theta_k \in (0,\pi)$, we have
$
\theta_k = \arccos \left(\frac{\cos (\xi_k - \xi_{k+1})+1}{2}\right).
$
Let 
$$
q_k := ( \cos \xi_k, \sin \xi_k, 1) \in \gamma_0.
$$
Then $q_k \to q= (1,0,1) \in \gamma_0$.

We now define $t_k$ depending on a parameter $a \in (-\frac{\pi}{2},\frac{\pi}{2}]$ by 
\begin{equation}\label{eq:tk_a}
t_k(a) = \frac{1}{\cos \left(a - \sum_{i=k}^\infty \theta_i \right)}, \quad k \geq 1.
\end{equation}
The limit behavior of points $p_k$ are given by the following lemma.
\begin{lemma}\label{lem:angle_conv_c}
  \begin{itemize}
      \item[1)] 
      For any $a \in \left(-\frac{\pi}{2}, \frac{\pi}{2}\right]$, there exists an integer $k_0(a)$ such that $t_k(a)$ is well-defined 
      (i.e., $\cos \left(a - \sum_{i=k}^\infty \theta_i \right) \neq 0$)
      and positive for all $k \geq k_0(a)$.
      \item[2)] 
      If $a \in \left(-\frac{\pi}{2}, \frac{\pi}{2}\right)$, then $t_k(a)$ converges to $\frac{1}{\cos a}$ as $k \to \infty$. Consequently, $p_k$ converges to a point $p=(\frac{1}{\cos a}, 0, \frac{1}{\cos a}) \in K_0$.
      \item[3)] 
      If $a = \frac{\pi}{2}$, then $t_k(a)$ tends to infinity as $k \to \infty$.
  \end{itemize}
\end{lemma}
\textbf{Proof.}
1)
Using \eqref{eq:xi-theta-0} and the half-angle formula we obtain
$$
\cos \theta_k
=
1 - 2 \sin^2 \frac{\theta_k}{2} 
= 1 - \sin^2 \frac{\xi_k -\xi_{k+1}}{2},
$$
hence
$$
\sin^2 \frac{\theta_k}{2} = \frac{1}{2}\sin^2 \frac{\xi_k -\xi_{k+1}}{2}.
$$
Since
$\frac{\theta_k}{2}, \frac{\xi_k - \xi_{k+1}}{2}\in (0,\frac{\pi}{2})$,  we have
\begin{equation}\label{eq:xi-theta-1}
\sin \frac{\theta_k}{2}= \frac{\sqrt{2}}{2} \sin \frac{\xi_k - \xi_{k+1}}{2}.
\end{equation}
By the monotonicity of $\sin x$ on $(0,\frac{\pi}{2})$, it follows that
$$
\frac{\theta_k}{2} < \frac{\xi_k - \xi_{k+1}}{2}.
$$

Therefore, for any $k,k'>0$,
$$
\sum_{i=k}^{k+k'} \theta_i < \sum_{i=k}^{k+k'} (\xi_i - \xi_{i+1}) = \xi_k - \xi_{k+k'+1} < k^{-1/2}.
$$
Letting $k'\to\infty$, we obtain that $\sum_{i=k}^\infty \theta_i$ converges and its sum does not exceed $k^{-1/2}$. 
In particular,
\begin{equation}\label{eq:theta-sum-limit}
\lim_{k\to\infty}\sum_{i=k}^\infty \theta_i = 0.
\end{equation}

Now, for any $a \in \left(-\frac{\pi}{2}, \frac{\pi}{2}\right]$, since $\sum_{i=k}^\infty \theta_i < k^{-1/2}$, we can choose an integer $k_0(a)$ such that 
$$
a - k_0(a)^{-1/2} >  -\frac{\pi}{2}.
$$ 
It then follows that for all $k \geq k_0(a)$, 
$$
a-\sum_{i=k}^\infty \theta_i \in (-\frac{\pi}{2}, \frac{\pi}{2}).
$$
Thus, $t_k(a)$ is well-defined and positive for all $k \geq k_0(a)$ by (\ref{eq:tk_a}).

2) 
If $a \in \left(-\frac{\pi}{2}, \frac{\pi}{2}\right)$, it follows from \eqref{eq:tk_a} and \eqref{eq:theta-sum-limit} that
$$
\lim_{k \to \infty} t_k(a) = \frac{1}{\cos a}.
$$
Together with $\lim_{k\to\infty}\xi_k = 0$, this implies that $p_k$ converges to 
$p=(\frac{1}{\cos a}, 0, \frac{1}{\cos a})\in K_0$.

3)
If $a = \frac{\pi}{2}$, we have
$$
t_k\left(\frac{\pi}{2}\right) = \frac{1}{\cos\left(\frac{\pi}{2} - \sum_{i=k}^\infty \theta_i\right)} \to \infty 
\quad \text{as}\ k \to \infty.
$$

Lemma \ref{lem:angle_conv_c} is proved.

\vspace*{1em}

The sequence $p_k, k\geq k_0(a)$ defines a polygonal line inside the cone $K_0$. Let $l_k$ denote the oriented line that contains the segment from $p_k$ to $p_{k+1}$, 
and let $v_{k} = \frac{p_{k+1}-p_{k}}{\|p_{k+1}-p_{k}\|}$ be its direction.
We have the following.

\begin{lemma}\label{lem:trajectory_property}
  \begin{itemize}
    \item[1)] The distance from $l_k$ to the origin equals $\sqrt{2}$ for all $k\geq k_0(a)$.
    \item[2)] The lines $l_{k-1}$ and $l_{k}$ form the same angle with the vector $p_{k}$ for all $k > k_0(a)$.
  \end{itemize}
\end{lemma}
\textbf{Proof.}
1)
Let us compute $\text{dist}(l_k, O)$. 
From (\ref{eqdistance}) we have 
$$
\begin{aligned}
\text{dist}(l_k, O)^2 
&= \|p_{k+1}\|^2 - \left\langle p_{k+1},  \frac{p_{k+1}-p_{k}}{\|p_{k+1}-p_{k}\|} \right \rangle^2\\
&= \|p_{k+1}\|^2 - 
\frac{
  \left(
  \|p_{k+1}\|^2 - \langle p_{k+1}, p_k \rangle
  \right)^2
  }
  {\|p_{k+1}\|^2 + \|p_k\|^2 - 2\langle p_{k+1}, p_k \rangle}\\
&= \|p_{k+1}\|^2
\left(1 - \frac{ \left(
  \|p_{k+1}\| - \| p_k\| \cos \theta_k 
  \right)^2}
  {\|p_{k+1}\|^2 + \| p_k\|^2 - 2\| p_{k+1}\| \|p_{k}\| \cos \theta_k}
\right)\\
&= \|p_{k+1}\|^2
\left( \frac{ 
   \| p_k\|^2 \sin^2 \theta_k 
  }
  {\|p_{k+1}\|^2 + \| p_k\|^2 - 2\| p_{k+1}\| \|p_{k}\| \cos \theta_k}
\right)\\
&= \frac{\sin^2 \theta_k }{\|p_{k+1}\|^{-2} + \| p_k\|^{-2} -2 \| p_{k+1}\|^{-1} \|p_{k}\|^{-1} \cos \theta_k}
.
\end{aligned}
$$
From (\ref{eq:def-of-pk}), (\ref{eq:tk_a}), we have 
\begin{equation}\label{eq:distance-pk}
  \|p_k\|^{-1}
  =\left(\sqrt{2} t_k(a)\right)^{-1}
  =
 \frac{\cos\left(a - \sum_{i=k}^\infty \theta_i\right)}{\sqrt{2}},
  \quad 
\forall k\geq k_0(a).
\end{equation}
Substituting (\ref{eq:distance-pk}) into $2\left(\|p_{k+1}\|^{-2} + \| p_k\|^{-2} -2 \| p_{k+1}\|^{-1} \|p_{k}\|^{-1} \cos \theta_k\right)$ we obtain
$$
\begin{aligned}
&
2\left(\|p_{k+1}\|^{-2} + \| p_{k}\|^{-2} -2 \| p_{k+1}\|^{-1} \|p_{k}\|^{-1} \cos \theta_k\right)
\\
=\ & \cos^2 \left(a- \sum_{i={k+1}}^{\infty} \theta_i \right) 
+\cos^2 \left(a- \sum_{i=k}^{\infty} \theta_i\right) 
- 2\cos \left(a-  \sum_{i=k+1}^{\infty} \theta_i\right) 
\cos \left(a-  \sum_{i=k}^{\infty} \theta_i\right)  \cos \theta_k
\\
=\ & \cos \left(a- \sum_{i={k+1}}^{\infty} \theta_i \right) 
\cos \left(\left(a- \sum_{i={k}}^{\infty} \theta_i \right)+ \theta_k \right) 
+ \cos \left(a- \sum_{i=k}^{\infty} \theta_i\right) 
\cos \left(\left(a- \sum_{i=k+1}^{\infty} \theta_i \right)- \theta_k\right) \\
&- 2\cos \left(a-  \sum_{i=k+1}^{\infty} \theta_i\right) 
\cos \left(a-  \sum_{i=k}^{\infty} \theta_i\right)  \cos \theta_k\\
=\ &\cos \left(a- \sum_{i={k+1}}^{\infty} \theta_i \right) 
\left(\cos 
\left(a- \sum_{i={k}}^{\infty} \theta_i \right) 
\cos \theta_k 
- \sin
\left(a- \sum_{i={k}}^{\infty} \theta_i \right)  
 \sin \theta_k \right)
 \\
& +
\cos 
\left(a- \sum_{i={k}}^{\infty} \theta_i \right) 
\left(\cos
\left(a- \sum_{i={k+1}}^{\infty} \theta_i \right) 
 \cos \theta_k 
+ \sin
\left(a- \sum_{i={k+1}}^{\infty} \theta_i \right) 
\sin \theta_k \right)\\
& - 2\cos \left(a-  \sum_{i=k+1}^{\infty} \theta_i\right) 
\cos \left(a-  \sum_{i=k}^{\infty} \theta_i\right)  \cos \theta_k
\\
=\ &\left(
\sin 
\left(a- \sum_{i={k+1}}^{\infty} \theta_i \right) 
\cos 
\left(a- \sum_{i={k}}^{\infty} \theta_i \right) 
-
\cos 
\left(a- \sum_{i={k+1}}^{\infty} \theta_i \right) 
\sin
\left(a- \sum_{i={k}}^{\infty} \theta_i \right)
\right)
\sin \theta_k\\
=\ &\sin^2\theta_k .
\end{aligned}
$$
Therefore, 
$$
\text{dist}(l_k, O)^2 = 2, \quad \forall k\geq k_0(a).
$$

\begin{figure}[htbp]
  \begin{center}
  \includegraphics[scale=0.48]{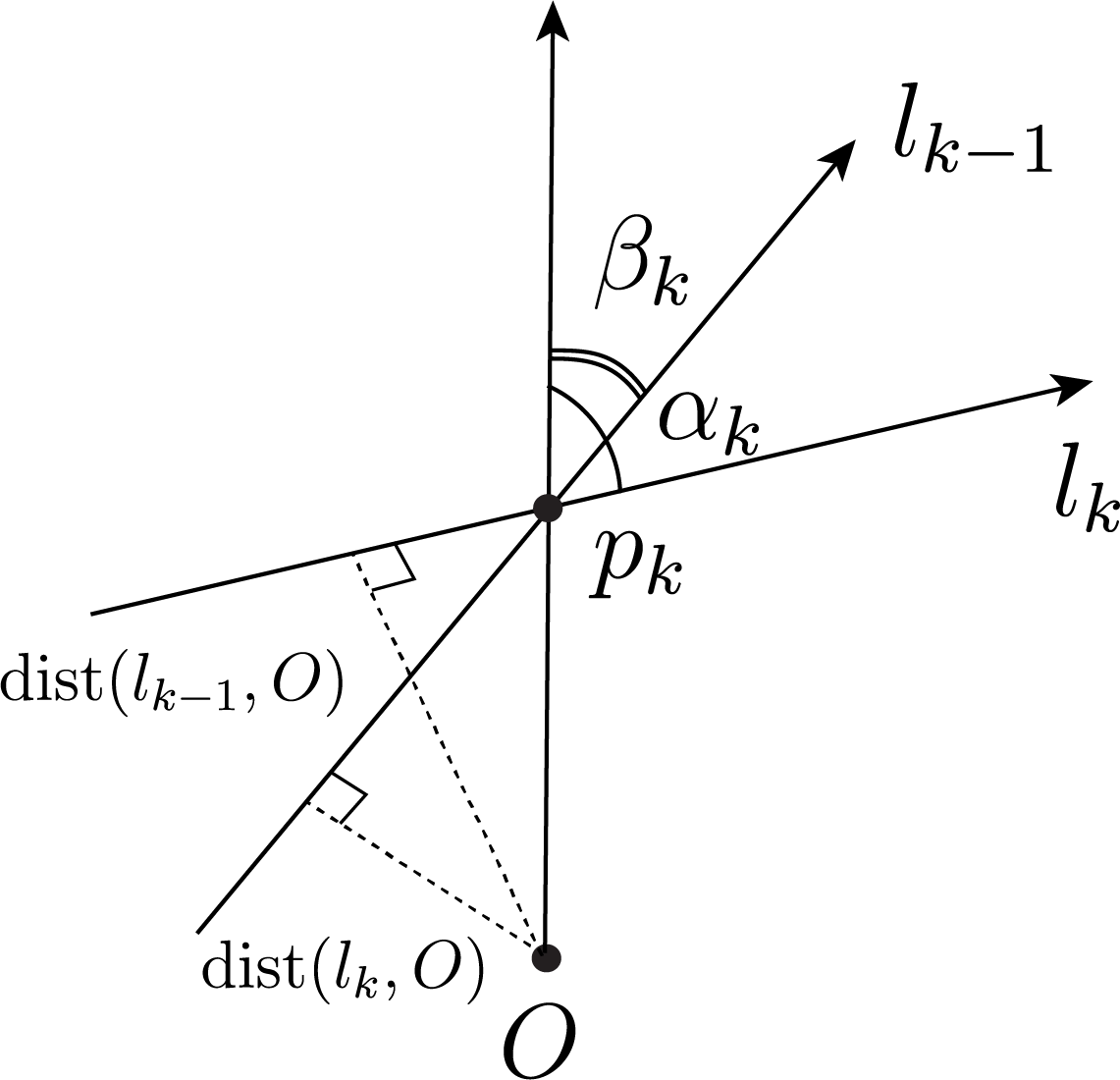}
  \end{center}
  \captionsetup{width=0.75\textwidth} 
  \caption{The angle $\beta_{k}$ between $l_{k-1}$ and $Op_{k}$. (Lines $l_k$, $l_{k-1}$ and $Op_{k}$ do not lie in the same plane.) }
  \label{fig:alpha-beta}
\end{figure}

2)
Let $\beta_{k}$ be the angle between $l_{k-1}$ and the vector $p_{k}$ (see Fig. \ref{fig:alpha-beta}). 
We will show $\cos \alpha_k= \cos \beta_k$.

The calculation of the area of $\triangle Op_k p_{k+1}$ in two ways gives
\begin{equation}\label{eq:area-p}
\frac{1}{2}\|p_{k+1} - p_{k}\| \text{dist}(l_{k}, O) =\frac{1}{2} \|p_{k+1}\|\|p_{k}\|\sin \theta_{k}.
\end{equation}
From part 1) we have $\text{dist}(l_{k}, O)=\sqrt{2}$, 
thus (\ref{eq:area-p}) gives
\begin{equation}\label{eq:lk-distance}
  \|p_{k+1} - p_{k}\| = \frac{\sin \theta_k}{\sqrt{2}}\|p_{k+1}\|\|p_{k}\|,
\quad 
k\geq k_0(a).
\end{equation}

Now let us compute $\cos \alpha_k$.
Using (\ref{eq:lk-distance}), we have
$$
\begin{aligned}
  \cos \alpha_k 
& = \left\langle \frac{p_{k+1} - p_k}{\|p_{k+1} - p_k\|}, \frac{p_k}{\|p_k\|} \right\rangle 
= \frac{\|p_{k+1}\|\|p_k\| \cos \theta_k - \|p_k\|^2}{\|p_{k+1} - p_k\|\|p_k\|}
= \frac{\|p_{k+1}\| \cos \theta_k - \|p_k\|}{\|p_{k+1} - p_k\|}\\
&= \frac{\sqrt{2}\ }{\sin \theta_k} \frac{\|p_{k+1}\| \cos \theta_k - \|p_k\|}{\|p_{k+1}\| \|p_k\| } 
= \frac{\sqrt{2}\ }{\sin \theta_k} 
\left(
  \|p_{k}\|^{-1} \cos \theta_k - \|p_{k+1}\|^{-1}
\right).
\end{aligned}
$$
Substituting (\ref{eq:distance-pk}) we obtain
\begin{equation}\label{eq:cos-alphak}
\begin{aligned}
\cos \alpha_k 
&= \frac{1}{\sin \theta_k} 
\left(
  \cos\left(a- \sum_{i=k}^\infty \theta_i\right)\cos \theta_k - \cos\left(a- \sum_{i=k+1}^\infty \theta_i\right)
\right)\\
&= 
\frac{1}{\sin \theta_k} 
\left(
  \cos\left(
    a- \sum_{i=k}^\infty \theta_i
    \right)\cos \theta_k 
    - \cos\left(a- \sum_{i=k}^\infty \theta_i
    +\theta_k
    \right)
\right)\\
&= 
\frac{1}{\sin \theta_k} 
\left(
  \cos\left(
    a- \sum_{i=k}^\infty \theta_i \right)
   \cos \theta_k 
   -
   \cos\left(a- \sum_{i=k}^\infty \theta_i\right)
   \cos \theta_k 
   +
   \sin\left(
    a- \sum_{i=k}^\infty \theta_i \right) 
    \sin \theta_k
\right)\\
&= \sin\left(
  a- \sum_{i=k}^\infty \theta_i \right) 
.
\end{aligned}
\end{equation}
Similarly, for $\cos \beta_k$ we have
$$
\cos \beta_k 
= \left\langle \frac{p_{k} - p_{k-1}}{\|p_{k} - p_{k-1}\|}, \frac{p_{k}}{\|p_{k}\|} \right\rangle 
= \frac{\|p_{k}\| - \|p_{k-1}\| \cos \theta_{k-1}}{\|p_{k} - p_{k-1}\|}
=\frac{\sqrt{2}\ }{\sin \theta_{k-1}} 
\left(
  \|p_{k-1}\|^{-1} - \|p_{k}\|^{-1}\cos \theta_{k-1}
\right),
$$
and then 
\begin{equation}
\label{eq:cos-betak}
\begin{aligned}
\cos \beta_k 
&= 
\frac{1}{\sin \theta_{k-1}} 
\left(
  \cos\left(
    a- \sum_{i=k}^\infty \theta_i - \theta_{k-1}
    \right)
    - \cos\left(a- \sum_{i=k}^\infty \theta_i\right)
    \cos \theta_{k-1} 
\right)\\
&= 
\frac{1}{\sin \theta_{k-1}} 
\left(
  \cos\left(
    a- \sum_{i=k}^\infty \theta_i \right)
   \cos \theta_{k-1}
   + \sin\left(
    a- \sum_{i=k}^\infty \theta_i \right) 
    \sin \theta_{k-1}
    - 
    \cos\left(a- \sum_{i=k}^\infty \theta_i\right)
    \cos \theta_{k-1} 
\right)\\
&= \sin\left(
  a- \sum_{i=k}^\infty \theta_i \right) 
.
\end{aligned}
\end{equation}
By (\ref{eq:cos-alphak}) and (\ref{eq:cos-betak}) we have
$\cos \alpha_{k}=\cos \beta_{k}$.
Since $\alpha_k, \beta_k \in (0, \pi)$, 
it follows that
$$
\alpha_{k}=\beta_{k}.
$$

Lemma \ref{lem:trajectory_property} is proved.

\begin{remark}
  From \eqref{eq:cos-alphak}, we observe that $\cos \alpha_k \to \sin a$ as $k \to \infty$. Let $\alpha$ be the limit of $\alpha_k$. This implies  $\cos \alpha =  \sin a $.
  Since $a \in \left(-\frac{\pi}{2}, \frac{\pi}{2}\right]$ and $\alpha \in [0, \pi)$, it follows that $\alpha = \frac{\pi}{2} - a$.  
  Geometrically, this means that the oriented line $l_k$ converges to a limiting oriented line $l$ forming an angle $\alpha = \frac{\pi}{2} - a$ with the oriented line $Oq$. 
\end{remark}

By part 2) of Lemma \ref{lem:angle_conv_c}, for  $a \in \left(-\frac{\pi}{2}, \frac{\pi}{2}\right)$, $p_k$ converges to $p=(\frac{1}{\cos a}, 0, \frac{1}{\cos a})$.
The next lemma gives the total length of the polygonal line determined by $p_k$, $k \geq k_0(a)$.

\begin{lemma}\label{cor:leng}
For $a \in \left(-\frac{\pi}{2}, \frac{\pi}{2}\right)$, we have  
\begin{equation*}
\sum_{k=k_0(a)}^{\infty} \|p_{k+1} - p_{k}\| =  
\frac{\sqrt{2} \sin \left(\sum_{i=k_0(a)}^{\infty} \theta_i\right) }{\cos \left(a - \sum_{i=k_0(a)}^{\infty} \theta_i\right) \cos a}.
\end{equation*}
\end{lemma}
\textbf{Proof.}
We calculate the area of $\triangle p_{k}Op_{k+1}$, $k \geq k_0(a)$ in two ways:
$$
\frac{1}{2} \|p_{k}\| \|p_{k+1}\| \sin \theta_k
= \frac{1}{2} \|p_{k+1}-p_{k}\| \text{dist}(l_k,O),
$$
where $\theta_k =\angle p_{k}Op_{k+1}$.
From (\ref{eq:distance-pk}) and Lemma \ref{lem:trajectory_property}, we have
$$
\|p_{k}\| =   \frac{\sqrt{2}}{\cos \left(a - \sum_{i=k}^{\infty} \theta_i\right)},
\quad 
\text{dist}(l_k,O)= \sqrt{2}.
$$
Hence we have
\begin{equation*}
  \|p_{k+1}-p_{k}\| =\frac{\|p_{k}\| \|p_{k+1}\| \sin \theta_k}{ \text{dist}(l_k,O)}
  = 
  \frac{\sqrt{2} \sin \theta_k }{
    \cos \left(a - \sum_{i=k+1}^{\infty} \theta_i\right) 
    \cos \left(a - \sum_{i=k}^{\infty} \theta_i\right)}
\end{equation*}
$$
=  \frac{\sqrt{2} 
\sin \left(
  \left(a - \sum_{i=k+1}^{\infty} \theta_i\right) 
  -
  \left(a - \sum_{i=k}^{\infty} \theta_i\right)
\right)
}{
  \cos \left(a - \sum_{i=k+1}^{\infty} \theta_i\right) 
  \cos \left(a - \sum_{i=k}^{\infty} \theta_i\right)}
=
\frac{\sqrt{2} 
\sin 
  \left(a - \sum_{i=k+1}^{\infty} \theta_i\right) 
}{
  \cos \left(a - \sum_{i=k+1}^{\infty} \theta_i\right)}
  -
  \frac{\sqrt{2} 
\sin
  \left(a - \sum_{i=k}^{\infty} \theta_i\right) 
}{
  \cos 
  \left(a - \sum_{i=k}^{\infty} \theta_i\right)}.
$$
Therefore, for $k' > k_0(a)$
\begin{equation*}
\sum_{k=k_0(a)}^{k'} \|p_{k+1} - p_{k}\|
=
\frac{\sqrt{2} 
\sin 
  \left(a - \sum_{i=k'+1}^{\infty} \theta_i\right) 
}{
  \cos \left(a - \sum_{i=k'+1}^{\infty} \theta_i\right)}
  -
  \frac{\sqrt{2} 
\sin
  \left(a - \sum_{i=k_0(a)}^{\infty} \theta_i\right) 
}{
  \cos \left(a - \sum_{i=k_0(a)}^{\infty} \theta_i\right)}.
\end{equation*}
Letting $k' \to \infty$ and using (\ref{eq:theta-sum-limit}), we obtain
\begin{equation*}
  \sum_{k=k_0(a)}^{\infty} \|p_{k+1} - p_{k}\|
  =
  \frac{\sqrt{2} 
  \sin a
  }{
    \cos a}
    -
    \frac{\sqrt{2} 
  \sin 
    \left(a - \sum_{i=k}^{\infty} \theta_i\right) 
  }{
    \cos \left(a - \sum_{i=k}^{\infty} \theta_i\right)}
  =
  \frac{\sqrt{2} \sin \left(\sum_{i=k_0(a)}^{\infty} \theta_i\right) }{\cos \left(a - \sum_{i=k_0(a)}^{\infty} \theta_i\right) \cos a}.
  \end{equation*}

Lemma \ref{cor:leng} is proved.  

\vspace*{1em}

Next, we perturb $K_0$ to obtain a new cone $K$ such that $p_k \in K$ and the polygonal line determined by $p_k$, $k \geq k_0(a)$ becomes a billiard trajectory inside $K$.
By Lemma \ref{cor:leng}, for $a \in \left(-\frac{\pi}{2}, \frac{\pi}{2}\right)$, this billiard trajectory is of finite length, and hence has infinitely many reflections in finite time.

\begin{remark}
  By part 3) of Lemma \ref{lem:angle_conv_c},
  $t_k(a)\to \infty$ if $a=\frac{\pi}{2}$.
  In this case
  the polygonal line defined by $p_k$, $k\geq k_0(a)$ has infinite length. 
  Therefore, by the same construction, one can also prove that 
  there exist $C^2$-smooth convex cones with billiard trajectories having infinitely many reflections in infinite time.
\end{remark}

To construct $K$, we construct a strictly convex $C^2$ curve $\gamma \subset \mathcal{P}$ with $q_k = (\cos \xi_k, \sin \xi_k, 1) \in \gamma$. The cone $K$ will be the cone over $\gamma$.

Let $T_k$ be the unit tangent vector of $\gamma$ at $q_k$ (hence $T_k \in T_{p_k}K \subset \mathbb{R}^n$). 
The polygonal line determined by $p_k$, $k \geq k_0(a)$ becomes a billiard trajectory in $K$ if $T_k$ satisfies:
\begin{equation}\label{tk}
\langle T_k, {v}_{k} \rangle = \langle T_k, {v}_{k-1} \rangle.
\end{equation}
Indeed, from (\ref{tk}) and by
Lemma \ref{lem:trajectory_property}, part 2)
$$
\langle p_k, {v}_{k} \rangle = \langle p_k, {v}_{k-1} \rangle, 
\quad p_k \in T_{p_k} K.
$$
Hence, the angles between ${v}_{k-1}$ and $T_{p_k}K$ and between ${v}_{k}$ and $T_{p_k}K$ are equal.

From (\ref{tk}) we have
$$
\langle T_k, {v}_{k} - {v}_{k-1} \rangle = 0.
$$
Define
$$
w_k := \frac{\tilde{w}_k}{\|\tilde{w}_k\|},
$$
where $\tilde{w}_k$ is the projection of ${v}_{k} - {v}_{k-1}$ onto $\mathcal{P}$. The vector $w_k$ serves as the normal vector to $\gamma$ at $q_k$. We then take $T_k$ to be the vector obtained by rotating $w_k$ clockwise by $\pi/2$ in the plane $\mathcal{P}$. 
Then $T_k$ satisfies condition (\ref{tk}).

Finally, we need to construct a $C^2$-smooth strictly convex plannar curve $\gamma$ 
that passes through all the points $q_k$ and has $w_k$ as the inward normal vector at each $q_k$. 
The existence of such a curve $\gamma$ is guaranteed by the following Lemmas \ref{lem:angle_convergence} and \ref{lem:Hal}.

\begin{figure}[htbp]
  \begin{center}
  \includegraphics[scale=0.28]{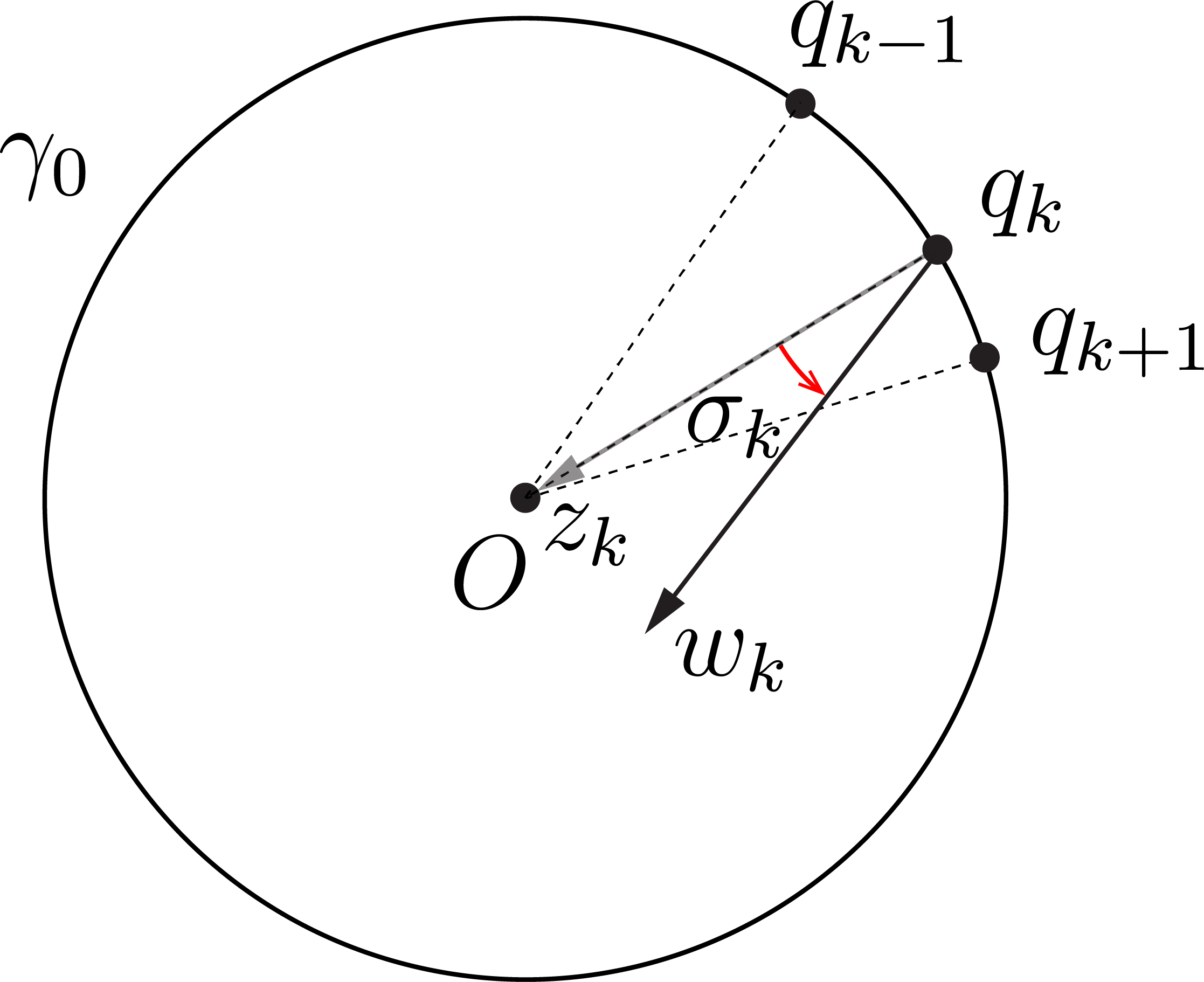}
  \end{center}
  \caption{The oriented angle $\sigma_k$ measured counterclockwise from $ z_k $ to $ w_k $.}
  \label{fig:angle-sigma}
\end{figure}

Let $z_k$ be the unit inward normal vector to the unit curcle $\gamma_0\subset \mathcal{P}$ at $q_k$, 
and $\sigma_k$ be the oriented angle measured counterclockwise from $ z_k $ to $ w_k $
(see Fig. \ref{fig:angle-sigma}). 
\begin{lemma}\label{lem:angle_convergence}
We have the identity
  $$
  \sigma_k = b_k k^{-5/2},
  $$
  where $b_k \to \frac{3}{16}$ as $k \to \infty$.
\end{lemma}
\textbf{Proof.}
In the plane $\mathcal{P}$, the inward normal vector of $\gamma_0$ at $q_k$ is given by
$$
z_k = (-\cos \xi_k, -\sin \xi_k).
$$
Let $\tilde{z}_k$ be the vector obtained by rotating $z_k$ counterclockwise by $\frac{\pi}{2}$ in $\mathcal{P}$, i.e.,
$$
\tilde{z}_k = (\sin \xi_k,-\cos \xi_k).
$$
Since $\sigma_k$ is the oriented angle measured counterclockwise from $z_k$ to $w_k$, 
$w_k$ can be expressed as
\begin{equation}\label{eq:wk_decomp}
w_k = \cos\sigma_k z_k + \sin\sigma_k \tilde{z}_k.
\end{equation}

From (\ref{eq:sin-p1}) and $\text{dist}(l_k,O)= \sqrt{2}$, we have
\begin{equation}\label{eq:pk-sin-alpha}
  \|p_k\|= \frac{\sqrt{2}}{\sin \alpha_k},
\end{equation}
Using  (\ref{eq:lk-distance}), (\ref{eq:pk-sin-alpha}),
$\alpha_{k+1} = \alpha_{k} - \theta_{k}$, and
$\frac{p_k}{\|p_k\|}= \frac{q_k}{\sqrt{2}}$ (since $q_k\in \gamma_0$)
we obtain 
\begin{equation*}
  \begin{split}
    {v}_{k}-{v}_{k-1} &= \frac{p_{k+1} - p_k}{\|p_{k+1} - p_k\|} - \frac{p_k - p_{k-1}}{\|p_k - p_{k-1}\|} \\
   & = 
   \frac{\sin \alpha_k \sin \alpha_{k+1 }\left(p_{k+1} - p_k\right)}{\sqrt{2} \sin {\theta_k}}
   -
   \frac{\sin \alpha_{k-1} \sin \alpha_{k }\left(p_{k} - p_{k-1}\right)}{\sqrt{2} \sin {\theta_{k-1}}}\\
   & = 
\frac{\sin \alpha_{k}}{\sin {\theta_{k}}} \frac{p_{k+1}}{\|p_{k+1}\|}
- 
\frac{\sin \alpha_{k+1}}{\sin {\theta_{k}}} \frac{p_{k}}{\|p_{k}\|}
-\frac{\sin \alpha_{k-1}}{\sin {\theta_{k-1}}} \frac{p_k}{\|p_{k}\|}
+ \frac{\sin \alpha_{k}}{\sin {\theta_{k-1}}} \frac{p_{k-1}}{\|p_{k-1}\|}\\
& = 
\frac{\sin \alpha_{k}}{\sin {\theta_{k}}} \frac{q_{k+1}}{\sqrt{2}}
- 
\frac{\sin \alpha_{k+1}}{\sin {\theta_{k}}} \frac{q_{k}}{\sqrt{2}}
-\frac{\sin \alpha_{k-1}}{\sin {\theta_{k-1}}} \frac{q_k}{\sqrt{2}}
+ \frac{\sin \alpha_{k}}{\sin {\theta_{k-1}}} \frac{q_{k-1}}{\sqrt{2}}\\
&= 
\frac{\sin \alpha_{k}}{\sin {\theta_{k}}} \frac{q_{k+1}}{\sqrt{2}}
+  \frac{\sin \alpha_{k}}{\sin {\theta_{k-1}}} \frac{q_{k-1}}{\sqrt{2}}
- \left(
  \frac{\sin (\alpha_{k}-\theta_k)}{\sin {\theta_{k}}}
+  \frac{\sin (\alpha_{k}+ \theta_{k-1})}{\sin {\theta_{k-1}}}
\right)
\frac{q_{k}}{\sqrt{2}}\\
& =\frac{\sin \alpha_{k}}{\sin {\theta_{k}}} \frac{q_{k+1}}{\sqrt{2}}
+  \frac{\sin \alpha_{k}}{\sin {\theta_{k-1}}} \frac{q_{k-1}}{\sqrt{2}}
- \left(
  \frac{\sin \alpha_{k} \cos \theta_k}{\sin {\theta_{k}}}
+  \frac{\sin \alpha_{k} \cos \theta_{k-1} }{\sin {\theta_{k-1}}}
\right)
\frac{q_{k}}{\sqrt{2}}.
  \end{split}
\end{equation*}
To compute the projection of $ {v}_{k}-{v}_{k-1}$ to $\mathcal{P}$, we introduce $N_k$ which differs from $ {v}_{k}-{v}_{k-1}$ only by a positive scalar factor. This simplifies the computation:
\begin{equation}\label{eq:nk}
  N_k=(N_{k}^1,N_{k}^2,N_{k}^3):= \frac{\sqrt{2} }{\sin \alpha_k}
   (  {v}_{k}-{v}_{k-1})
=
\frac{q_{k+1}}{\sin {\theta_{k}}} 
+  \frac{q_{k-1}}{\sin {\theta_{k-1}}} 
- \left(
  \frac{\cos \theta_k}{\sin {\theta_{k}}}
+  \frac{\cos \theta_{k-1} }{\sin {\theta_{k-1}}}
\right)
q_{k} .
\end{equation}
Substituting coordinates of $q_{k-1}, q_k, q_{k+1}$, 
and letting
$\delta_{k}:= \xi_{k}-\xi_{k+1}$
\begin{equation*}
  \begin{split}
N_{k}^1 
& = 
\frac{\cos \xi_{k+1}}{\sin {\theta_{k}}} 
+  \frac{\cos \xi_{k-1}}{\sin {\theta_{k-1}}} 
- \left(
  \frac{\cos \theta_k}{\sin {\theta_{k}}}
+  \frac{\cos \theta_{k-1} }{\sin {\theta_{k-1}}} 
\right)
\cos \xi_{k}\\
&= \frac{\cos (\xi_{k}-\delta_k)}{\sin {\theta_{k}}} 
+  \frac{\cos (\xi_{k}+\delta_{k-1})}{\sin {\theta_{k-1}}} 
- \left(
  \frac{\cos \theta_k}{\sin {\theta_{k}}}
+  \frac{\cos \theta_{k-1} }{\sin {\theta_{k-1}}} 
\right)
\cos \xi_{k}\\
&= \frac{\cos \xi_{k} \cos \delta_k}{\sin {\theta_{k}}} 
+\frac{\sin\xi_{k} \sin\delta_k}{\sin {\theta_{k}}} 
+  \frac{\cos\xi_{k} \cos \delta_{k-1}}{\sin {\theta_{k-1}}} 
-  \frac{\sin\xi_{k} \sin\delta_{k-1}}{\sin {\theta_{k-1}}} 
- \left(
  \frac{\cos \theta_k}{\sin {\theta_{k}}}
+  \frac{\cos \theta_{k-1} }{\sin {\theta_{k-1}}} 
\right)
\cos \xi_{k}\\
&= \left(
\frac{\cos \delta_{k}-\cos \theta_k}{\sin \theta_k}
+ \frac{\cos \delta_{k-1}-\cos \theta_{k-1}}{\sin \theta_{k-1}}
\right) \cos \xi_{k}
+ \left(
  \frac{\sin \delta_{k}}{\sin \theta_k}
- \frac{\sin \delta_{k-1}}{\sin \theta_{k-1}}
\right) \sin \xi_k .
\end{split}
\end{equation*}
From (\ref{eq:xi-theta-1}) we have 
\begin{equation}\label{eq:th-del}
  \sin \frac{\theta_k}{2}= \frac{\sqrt{2}}{2} \sin \frac{\delta_k}{2}.
\end{equation}
Then using the half-angle formula and (\ref{eq:th-del}) we simplify $N_k^1$ as
\begin{equation*}
  \begin{split}
N_{k}^1 
= \ & \left(
  \frac{\left(1- 2 \sin^2 \frac{\delta_{k}}{2}\right)-\left(1 -2\sin^2 \frac{\theta_k}{2} \right)}{2 \sin \frac{\theta_k}{2} \cos\frac{\theta_k}{2}}
  + 
  \frac{\left(1- 2 \sin^2 \frac{\delta_{k-1}}{2}\right)-\left(1 -2\sin^2 \frac{\theta_{k-1}}{2} \right)}{2 \sin \frac{\theta_{k-1}}{2} \cos\frac{\theta_{k-1}}{2}}
  \right) \cos \xi_{k}\\
  & + 
  \left(
    \frac{2\sin \frac{\delta_{k}}{2} \cos\frac{\delta_{k}}{2}}{2\sin \frac{\theta_k}{2} \cos\frac{\theta_k}{2}}
  - 
  \frac{2\sin \frac{\delta_{k-1}}{2} \cos\frac{\delta_{k-1}}{2}}{2\sin \frac{\theta_{k-1}}{2} \cos\frac{\theta_{k-1}}{2}}
  \right) 
  \sin \xi_k \\
= \ & \left(
  \frac{-  \sin^2 \frac{\delta_{k}}{2}
  + \sin^2 \frac{\theta_k}{2}}{ \sin \frac{\theta_k}{2} \cos\frac{\theta_k}{2}}
  + 
  \frac{-  \sin^2 \frac{\delta_{k-1}}{2}
  + \sin^2 \frac{\theta_{k-1}}{2} }{ \sin \frac{\theta_{k-1}}{2} \cos\frac{\theta_{k-1}}{2}}
  \right) \cos \xi_{k} 
  + 
  \left(
    \frac{\sin \frac{\delta_{k}}{2} \cos\frac{\delta_{k}}{2}}{\sin \frac{\theta_k}{2} \cos\frac{\theta_k}{2}}
  - 
  \frac{\sin \frac{\delta_{k-1}}{2} \cos\frac{\delta_{k-1}}{2}}{\sin \frac{\theta_{k-1}}{2} \cos\frac{\theta_{k-1}}{2}}
  \right) 
  \sin \xi_k \\
= \ & \left(
  \frac{-  2\sin^2 \frac{\theta_{k}}{2}
  + \sin^2 \frac{\theta_k}{2}}{ \sin \frac{\theta_k}{2} \cos\frac{\theta_k}{2}}
  + 
  \frac{-  2\sin^2 \frac{\theta_{k-1}}{2}
  + \sin^2 \frac{\theta_{k-1}}{2} }{ \sin \frac{\theta_{k-1}}{2} \cos\frac{\theta_{k-1}}{2}}
  \right) \cos \xi_{k} \\
  & + 
  \left(
    \frac{\sqrt{2} \sin \frac{\theta_{k}}{2} \cos\frac{\delta_{k}}{2}}{\sin \frac{\theta_k}{2} \cos\frac{\theta_k}{2}}
  - 
  \frac{\sqrt{2}\sin \frac{\theta_{k-1}}{2} \cos\frac{\delta_{k-1}}{2}}{\sin \frac{\theta_{k-1}}{2} \cos\frac{\theta_{k-1}}{2}}
  \right) 
  \sin \xi_k \\
=\ & -\left(
  \frac{\sin  \frac{\theta_k}{2} }{\cos \frac{\theta_k}{2}}
  + \frac{\sin  \frac{\theta_{k-1}}{2}}{\cos \frac{\theta_{k-1}}{2}}
  \right) \cos \xi_{k}
  + \sqrt{2} \left(
    \frac{\cos \frac{\delta_{k}}{2}}{\cos \frac{\theta_k}{2}}
  - \frac{\cos \frac{\delta_{k-1}}{2}}{\cos \frac{\theta_{k-1}}{2}}
  \right) \sin \xi_k .
\end{split}
\end{equation*}
Similarily,
\begin{equation*}
  \begin{split}
N_{k}^2 
&= 
\frac{\sin \xi_{k+1}}{\sin {\theta_{k}}} 
+  \frac{\sin \xi_{k-1}}{\sin {\theta_{k-1}}} 
- \left(
  \frac{\cos \theta_k}{\sin {\theta_{k}}}
+  \frac{\cos \theta_{k-1} }{\sin {\theta_{k-1}}}
\right)
\sin \xi_{k} \\
& = 
\frac{\sin (\xi_{k}-\delta_k)}{\sin {\theta_{k}}} 
+  \frac{\sin (\xi_{k}+\delta_{k-1})}{\sin {\theta_{k-1}}} 
- \left(
  \frac{\cos \theta_k}{\sin {\theta_{k}}}
+  \frac{\cos \theta_{k-1} }{\sin {\theta_{k-1}}}
\right)
\sin \xi_{k} \\
& = 
\frac{\sin\xi_{k} \cos\delta_k}{\sin {\theta_{k}}} 
-\frac{\cos\xi_{k} \sin\delta_k}{\sin {\theta_{k}}} 
+  \frac{\sin \xi_{k}\cos\delta_{k-1}}{\sin {\theta_{k-1}}} 
+  \frac{\cos \xi_{k}\sin\delta_{k-1}}{\sin {\theta_{k-1}}}
- \left(
  \frac{\cos \theta_k}{\sin {\theta_{k}}}
+  \frac{\cos \theta_{k-1} }{\sin {\theta_{k-1}}}
\right)
\sin \xi_{k} \\
& = 
 \left(
\frac{\cos \delta_{k}-\cos \theta_k}{\sin \theta_k}
+ \frac{\cos \delta_{k-1}-\cos \theta_{k-1}}{\sin \theta_{k-1}}
\right) \sin \xi_k
- \left(
  \frac{\sin \delta_{k}}{\sin \theta_k}
- \frac{\sin \delta_{k-1}}{\sin \theta_{k-1}}
\right) \cos \xi_k ,
\end{split}
\end{equation*}
which simplifies to 
\begin{equation*}
  \begin{split}
N_{k}^2 
& = -\left(
  \frac{\sin  \frac{\theta_k}{2} }{\cos \frac{\theta_k}{2}}
  + \frac{\sin  \frac{\theta_{k-1}}{2}}{\cos \frac{\theta_{k-1}}{2}}
  \right) \sin \xi_k 
  - 
  \sqrt{2}
  \left(
    \frac{\cos \frac{\delta_{k}}{2}}{\cos \frac{\theta_k}{2}}
  - \frac{\cos \frac{\delta_{k-1}}{2}}{\cos \frac{\theta_{k-1}}{2}}
  \right) \cos \xi_k .
\end{split}
\end{equation*}
Thus 
\begin{equation}\label{eq:nk2}
  \cos \frac{\theta_k}{2}
\cos \frac{\theta_{k-1}}{2}
 (N_{k}^1,N_{k}^2) = g_k(-\cos \xi_k, -\sin \xi_k) + f_k  (\sin \xi_k, -\cos \xi_k)
 =  g_k z_k + f_k \tilde{z}_k,
\end{equation}
where 
\begin{equation}\label{eq:gk-fk}
  g_k = 
  {\sin  \frac{\theta_k}{2} }
  {\cos \frac{\theta_{k-1}}{2}}
  + {\cos \frac{\theta_k}{2}}{\sin  \frac{\theta_{k-1}}{2}}
  ,
  \quad
  f_{k} =\sqrt{2}\left(
  {\cos \frac{\delta_{k}}{2}}\cos \frac{\theta_{k-1}}{2}
  -  {\cos \frac{\theta_k}{2}}{\cos \frac{\delta_{k-1}}{2}}
  \right).
\end{equation}
Since $\frac{\sqrt{2}}{\sin \alpha_k}>0$ in (\ref{eq:nk}), and $\cos \frac{\theta_k}{2}\cos \frac{\theta_{k-1}}{2}>0$ in (\ref{eq:nk2}),
we obtain 
$$
w_k = \frac{g_k }{\sqrt{g_k^2 + f_k^2}} z_k
+ \frac{f_k }{\sqrt{g_k^2 + f_k^2}} \tilde{z}_k.
$$
Comparing this with (\ref{eq:wk_decomp}), we obtain 
\begin{equation}\label{eq:sin-g-f}
  \sin \sigma_k = \frac{f_k }{\sqrt{g_k^2 + f_k^2}}.
\end{equation}
Using $\xi_{k}=k^{-1/2}$, $\delta_{k}= \xi_{k}-\xi_{k+1}$, and (\ref{eq:th-del}) we have
\begin{equation}\label{eq:delta-theta}
  \frac{\delta_k}{2} = \frac{1}{2}\left(\frac{1}{\sqrt{k}}-  \frac{1}{\sqrt{k+1}}\right),
\quad
\frac{\theta_{k}}{2} =  \arcsin \left(\frac{1}{\sqrt{2}} \sin \left(\frac{1}{2}\left(\frac{1}{\sqrt{k}}-  \frac{1}{\sqrt{k+1}}\right)\right)\right).
\end{equation}
Substituting (\ref{eq:gk-fk}), (\ref{eq:delta-theta}) into (\ref{eq:sin-g-f}),
we get by explicit calculation
$$
\sin \sigma_k = \frac{3}{16} k^{-5/2} + O(k^{-9/2}).
$$
Hence 
$$
\sigma_k =  \frac{3}{16} k^{-5/2} + O(k^{-9/2}).
$$

Lemma \ref{lem:angle_convergence} is proved.

\vspace*{1em}

\begin{figure}[htbp]
	\centering
	\begin{minipage}[t]{0.3\linewidth}
		\centering
		\includegraphics[width=0.9\textwidth]{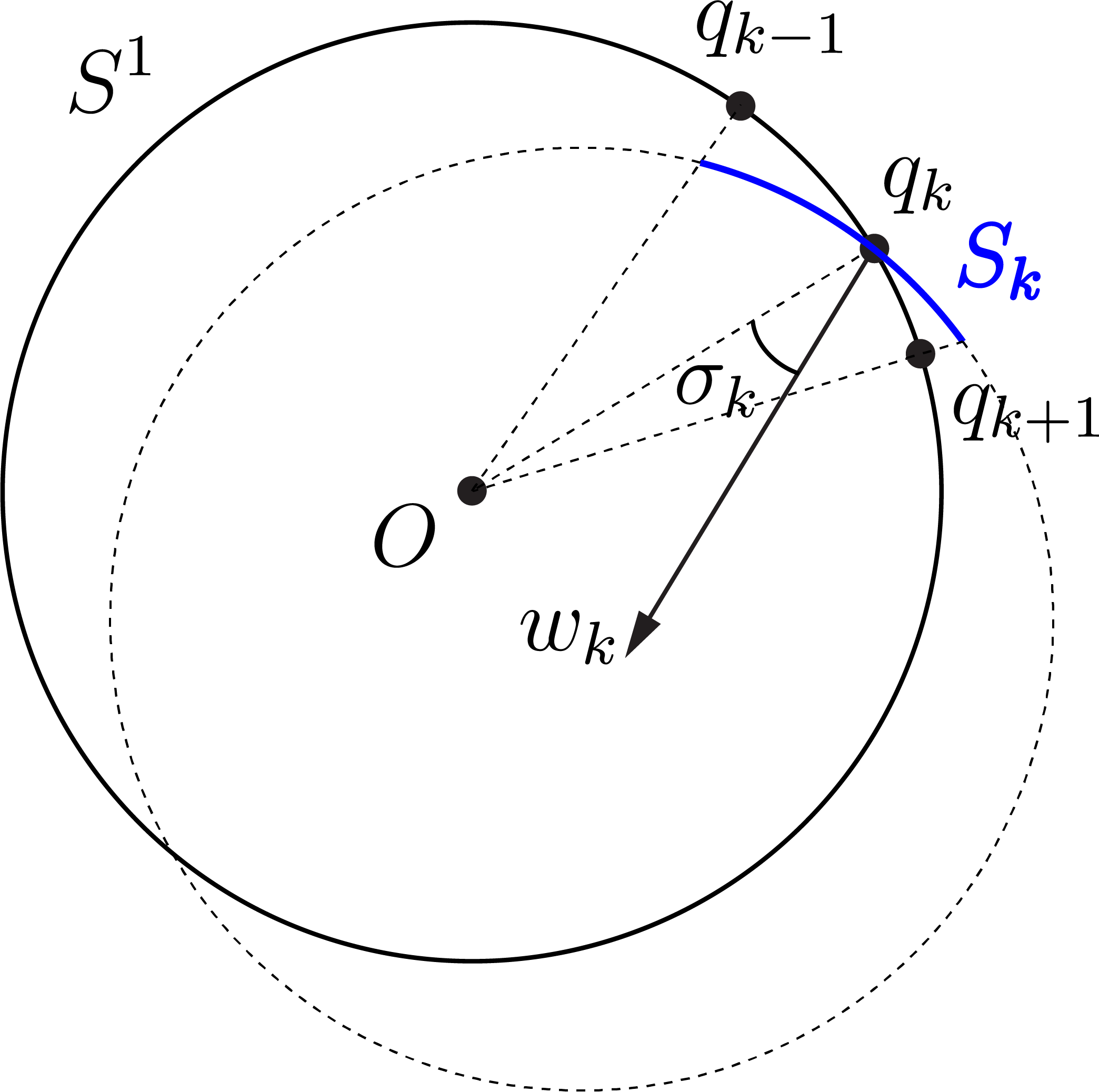}
    \caption{
      The arc $S_k$ near $q_k$.
    }
    \label{fig:sk}
	\end{minipage}
  \hspace{1cm}
  \begin{minipage}[t]{0.45\linewidth}
		\centering
		\includegraphics[width=0.9\textwidth]{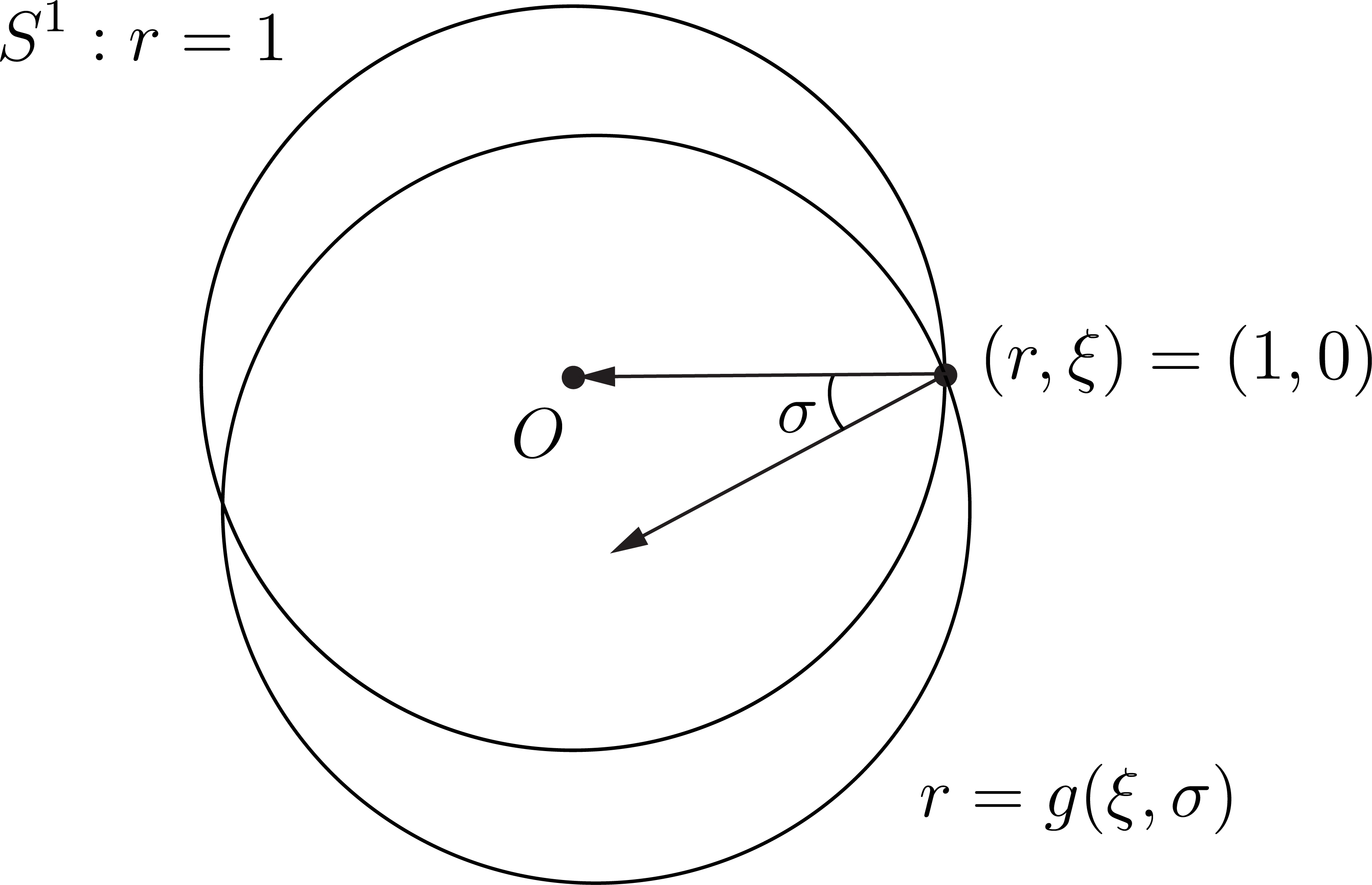}
    \caption{
The circle $r= g(\xi, \sigma)$. 
    }
    \label{fig:g-xi-sigma}
	\end{minipage}
\end{figure}

The following Lemma follows directly from Chapter 3 of \cite{Hal}, although it is not explicitly stated there. 
We give the proof of the lemma here for completeness.
\begin{lemma} {\bf (Halpern)}\label{lem:Hal}
  Consider a sequence of points 
  $q_k = (\cos \xi_k, \sin \xi_k)$ on the unit circle 
  $S^1 \subset \mathbb{R}^2$ with $\xi_k = k^{-1/2}$. 
  At each point $q_k$, a unit vector $w_k$ is assigned. 
  If the oriented angle $ \sigma_k $, measured counterclockwise from $z_k = -q_k $ to $ w_k $, satisfies  
  \begin{equation}\label{eq:asym}
      \left| \sigma_k \right| = b_k k^{-5/2}, \quad \text{with } b_k \to b > 0 \text{ as } k \to \infty,
  \end{equation}
  then there exists $k_0 > 0$ and a strictly convex $C^2$-smooth closed curve 
  $ \gamma \subset \mathbb{R}^2 $ such that $ \gamma $ passes through $q_k$ for all $k \geq k_0$, with the unit normal vector at each $q_k$ given by $w_k$.
\end{lemma}
\textbf{Proof.}
Let $(r, \xi)$ be the standard polar coordinates in $\mathbb{R}^2$, 
where $r$ is the radius and $\xi$ is the polar angle.
At each $q_k, k\geq 2$, draw a circle of radius $1$ passing through $q_k$ 
such that $w_k$ is the unit inward normal vector of the circle at $q_k$.
Consider an arc $S_k$ of this circle near $q_k$ (see Fig. \ref{fig:sk}), 
whose equation in the polar coordinates $(r,\xi)$ is given by 
$$
r = \rho_k(\xi),
\quad \xi \in \left[\xi_{k+1}, \xi_{k-1}\right], 
\quad  k\geq 2.
$$
Define $\rho_1(\xi)=1$, $\xi_2 \leq \xi \leq \xi_1$ 
(i.e., $S_1\subset S^1$). 
All the functions
$\rho_k(\xi), k\geq 1$ are $C^{\infty}$-smooth.
Let 
$$a(t): [0,1] \to [0,1]$$ 
be any $C^{\infty}$-smooth function such that $a(t) = 1$ for $t \leq \frac{1}{3}$, 
and $a(t) = 0$ for $t \geq \frac{2}{3}$.
Define $\rho(\xi)$ as follows:
\begin{equation}\label{eq:rho-rhok}
\rho(\xi) = \rho_k(\xi) + \left( \rho_{k+1}(\xi) - \rho_k(\xi) \right) a_k(\xi), 
\quad \text{for }\; \xi_{k+1} \leq \xi \leq \xi_k,
\end{equation}
where $a_k(\xi)$ is defined by
$$
a_k(\xi) = a\left( \frac{\xi - \xi_{k+1}}{\xi_k - \xi_{k+1}} \right).
$$
By this definition, we can see that $\rho(\xi)$ is well-defined on $(0,\xi_1]=(0,1]$, 
$C^{\infty}$-smooth on each $(\xi_{k+1},\xi_k)$.
$\rho(\xi)$ is also $C^{\infty}$-smooth at $\xi=\xi_k$, $k>1$, since 
\begin{equation}\label{eq:rhoatkp1}
\rho(\xi) = \rho_{k}(\xi), \quad \text{for } 
\xi \in 
\left[ \xi_{k}- \frac{\xi_k - \xi_{k+1}}{3}, \xi_{k} + \frac{\xi_{k-1} - \xi_{k}}{3} \right], 
\quad k > 1.
\end{equation}
Near $\xi=\xi_1=1$, we have
\begin{equation}
\label{eq:rhoatkp}
\rho(\xi) = 1, \quad \text{for } \xi \in 
\left[ \xi_1 - \frac{\xi_1 - \xi_2}{3}, \xi_1 \right].
\end{equation}
We extend $\rho(\xi)$ to $(-\pi, \pi]$ by setting 
\begin{equation}\label{eq:rhoatk1}
\rho(\xi) = 1, 
\quad
\text{for } \xi \in (-\pi, 0]\cup(\xi_1, \pi]. 
\end{equation}
It follows from (\ref{eq:rho-rhok})--(\ref{eq:rhoatk1})
that $\rho(\xi)$ is $C^{\infty}$-smooth on $(-\pi, 0)\cup(0,\pi]$.

Now, we define a curve $\tilde{\gamma}$ by
$$
r = \rho(\xi), \quad \xi \in (0, 2\pi].
$$
Then $\tilde{\gamma}$ is $C^{\infty}$-smooth on $(-\pi, 0)\cup(0,\pi]$.
We will prove that $\rho''(\xi)$ exists at $\xi = 0$ and is continuous, 
and thus $\tilde{\gamma}$ is $C^2$-smooth at $\xi = 0$.  
Moreover, we will prove that there exists a positive integer $k_0$ such that the 
curvature $\kappa(\xi) $ of  $\tilde{\gamma}$  satisfies
$$
\kappa(\xi) > \frac{1}{2},
\quad  \xi \in [0, \xi_{k_0}].
$$
We also mention that by the construction the vector $w_k$ is a 
normal vector to $\tilde{\gamma}$ at $q_k$.

\begin{figure}[htbp]
  \begin{center}
  \includegraphics[scale=0.3]{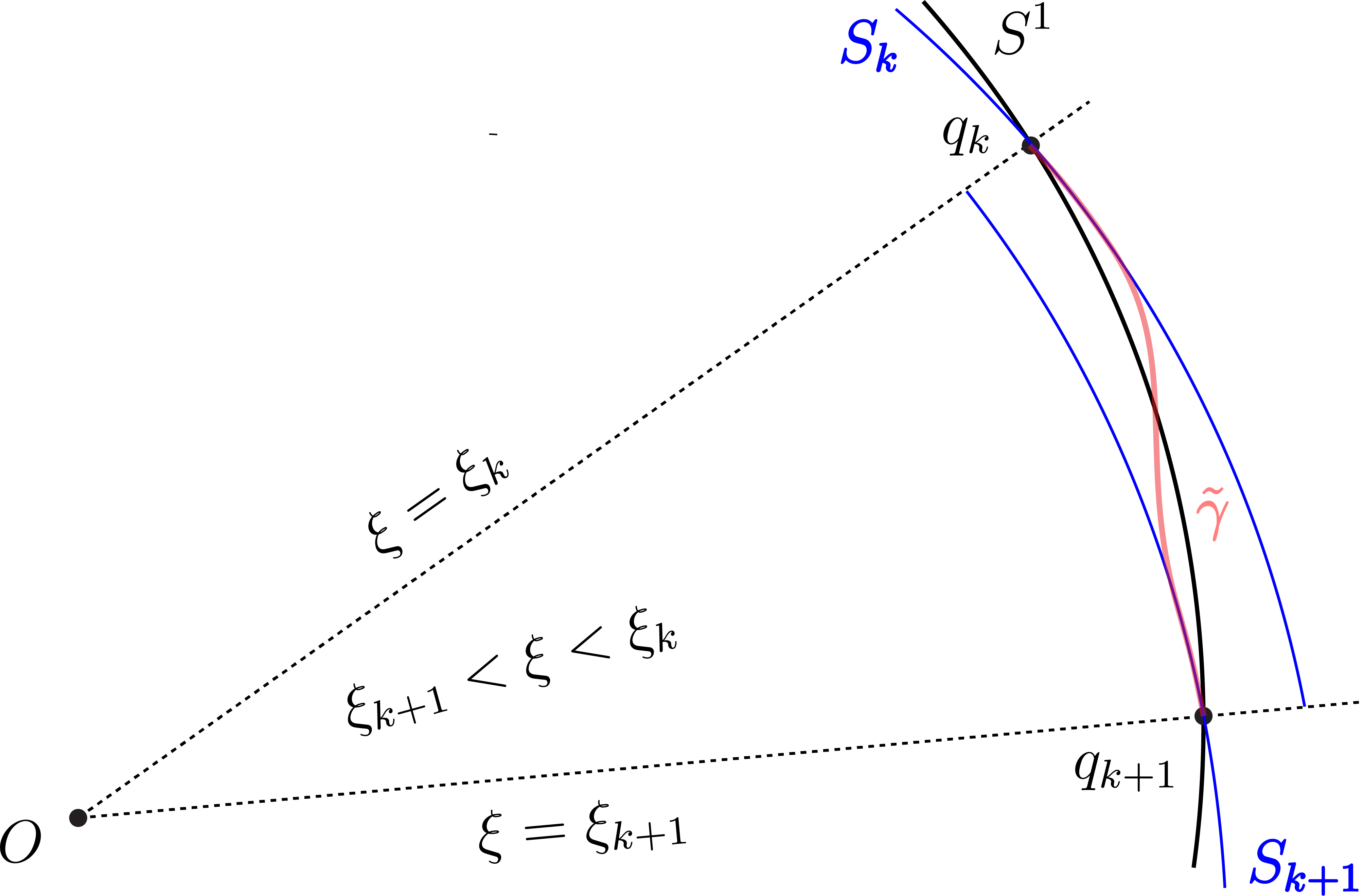}
  \end{center}
  \captionsetup{width=0.75\textwidth} 
  \caption{The blue curve passing through $q_k$ represents
   the arc $S_k$: $r = \rho_k(\xi)$, 
   $\xi_{k+1} \leq \xi \leq \xi_{k-1}$. 
  The blue curve passing through $q_{k+1}$ represents 
  the arc $S_{k+1}$: $r = \rho_{k+1}(\xi)$,  
  $\xi_{k+2} \leq \xi \leq \xi_k$. 
  The red curve passing through $q_k, q_{k+1}$ represents the curve $\tilde{\gamma}$: $r = \rho(\xi)$.
  $\tilde{\gamma}$
 coincides with $S_k$ near the point $q_k$, and coincides with $S_{k+1}$ near the point $q_{k+1}$.}
  \label{gamma}
\end{figure}

Using Taylor expansion, we obtain  
\begin{equation}\label{eq:delta-k}
  \delta_k := \xi_k - \xi_{k+1} = \frac{1}{2} k^{-3/2} + O(k^{-5/2}),  \text{ as } k \to \infty.
\end{equation}
Let $c>1$ be an upper bound for $|a'(t)|$ and $|a''(t)|$, $t \in [0,1]$. Then, for $k$ sufficiently large we have  
\begin{align}\label{eq:akp}
|a_k'(\xi)| 
=& \left|
\frac{a'\left(\frac{\xi - \xi_{k+1}}{\xi_k - \xi_{k+1}}\right)}{\xi_k - \xi_{k+1}} 
\right| 
= \frac{\left| a'\left(\frac{\xi - \xi_{k+1}}{\xi_k - \xi_{k+1}}\right)\right| }
{|\delta_k|}
\leq \frac{c}{|\delta_k|} \leq 3 c k^{3/2}, \\
\label{eq:akpp}
|a_k''(\xi)| 
=& \left|
\frac{a''\left(\frac{\xi - \xi_{k+1}}{\xi_k - \xi_{k+1}}\right)}{\left(\xi_k - \xi_{k+1}\right)^2} 
\right|
= \frac{\left| a''\left(\frac{\xi - \xi_{k+1}}{ \xi_k - \xi_{k+1}}\right)\right| }{\delta_k^2}
\leq \frac{c}{\delta_k^2} \leq 5 c k^{3}.
\end{align}

Let $r = g(\xi, \sigma)$ be the polar equation 
of the unit circle that passes through the point $(r, \xi) = (1, 0)$ and forms an angle $\sigma$ with $S^1$ at this point.
The angle $\sigma$ is measured counterclockwise from the inward normal of $S^1$ to the inward normal of the circle $r = g(\xi, \sigma)$ at $(r, \xi) = (1, 0)$ (see Fig. \ref{fig:g-xi-sigma}).
It is easy to verify that $g$ is well-defined and infinitely differentiable on the compact set
$$
\mathcal{D} = \{(\xi, \sigma) \mid -\pi/3 \leq \xi \leq \pi/3, -\pi/4 \leq \sigma \leq \pi/4\},
$$
and $g(0, \sigma) = g(\xi, 0)=1 $ for $(\xi, \sigma) \in \mathcal{D}$.

Let $e$ be a positive upper bound for 
$\left| \frac{\partial^2 g}{\partial \xi \partial \sigma} \right|$ and
$\left| \frac{\partial^3 g}{\partial \xi^2 \partial \sigma} \right|$ over domain $\mathcal{D}$. 
Then, by applying the mean value theorem, we obtain:  
\begin{equation}\label{eq:g-xi1}
\left| \frac{\partial g}{\partial \xi} (\xi, \sigma) - \frac{\partial g}{\partial \xi} (\xi, 0) \right| 
= \left| \frac{\partial^2 g}{\partial \xi \partial \sigma}(\xi,\sigma^*) (\sigma - 0) \right| 
= \left| \frac{\partial^2 g}{\partial \xi \partial \sigma}(\xi,\sigma^*)\right|\left|\sigma \right|,
\end{equation}
where $\sigma^*$ is some value between $0$ and $\sigma$.
Since $g(\xi, 0) = 1$, we have $\frac{\partial g}{\partial \xi} (\xi, 0) = 0$. 
Thus, (\ref{eq:g-xi1}) implies:
\begin{equation}\label{eq:g-xi2}
\left| \frac{\partial g}{\partial \xi}(\xi, \sigma) \right| \leq e \left| \sigma \right|.
\end{equation}
Similarly, we can derive:
\begin{equation}\label{eq:g-xi3}
\left| \frac{\partial^2 g}{\partial \xi^2 } (\xi, \sigma) \right| 
= \left| \frac{\partial^2 g}{\partial \xi^2} (\xi, \sigma) - \frac{\partial^2 g}{\partial \xi^2} (\xi, 0) \right| 
= \left| \frac{\partial^3 g}{\partial \xi^2 \partial \sigma}(\xi,\sigma^{**})\right|\left|\sigma \right| 
\leq e \left| \sigma \right|,
\end{equation}
where $\sigma^{**}$ is some value between $0$ and $\sigma$.
From $g(0, \sigma) = 1$ and  (\ref{eq:g-xi2}), we can derive: 
\begin{equation}\label{eq:g-xi4}
\left| g(\xi, \sigma) - 1 \right|  
= \left| g(\xi, \sigma) - g(0, \sigma) \right| 
= \left| \frac{\partial g}{\partial \xi} (\xi^*, \sigma) \right| \left| \xi \right| 
\leq e \left| \sigma \right| \left| \xi \right|,
\end{equation}
where $\xi^*$ is some value between $0$ and $\xi$.

By assumption (\ref{eq:asym}), 
there exists an integer $k_0\geq 2$ such that $|\sigma_k| < \frac{\pi}{4}$ for all $k \geq k_0$.
Since the circle $r= g(\xi - \xi_k, \sigma_k)$ is obtained by rotating the circle $r= g(\xi, \sigma_k)$ counterclockwise by angle $\xi_k$ about the origin, it forms an angle of $\sigma_k$ with $S^1$ at point $q_k$. Therefore, the arc $S_k: r= \rho_k(\xi), \xi \in \left[\xi_{k+1}, \xi_{k-1}\right]$ lies on the circle $r= g(\xi - \xi_k, \sigma_k)$ 
(compare Fig. \ref{fig:sk} and Fig. \ref{fig:g-xi-sigma}).
Thus, for $k \geq k_0$, we express $\rho_k(\xi)$ as
\begin{equation}\label{eq:rho-g} 
  \rho_k(\xi) = g(\xi - \xi_k, \sigma_k), 
  \quad 
  \xi \in \left[\xi_{k+1}, \xi_{k-1}\right]. 
\end{equation}
Here, since $|\xi - \xi_{k}|\leq |\xi_{k-1}- \xi_{k}|<1<\frac{\pi}{3}$ for $k\geq k_0$, we have $(\xi - \xi_k, \sigma_k)\in \mathcal{D}$.
Consequently, we obtain from  (\ref{eq:g-xi2}), (\ref{eq:asym}) that 
\begin{equation}\label{eq:rho-xi1}
\left| \rho_k'(\xi) \right| =
\left| \frac{\partial g}{\partial \xi}(\xi-\xi_k, \sigma_k) \right| 
 \leq e \left| \sigma_k \right| \leq 2 e b k^{-5/2},
\end{equation}
from  (\ref{eq:g-xi3}), (\ref{eq:asym}) that 
\begin{equation}\label{eq:rho-xi2}
\left| \rho_k''(\xi) \right| 
= \left| \frac{\partial^2 g}{\partial \xi^2 } (\xi-\xi_k, \sigma_k) \right| 
\leq e \left| \sigma_k \right| \leq 2 e b k^{-5/2},
\end{equation}
from (\ref{eq:g-xi4}), (\ref{eq:delta-k}), (\ref{eq:asym})
that
\begin{equation}\label{eq:rho-xi3}
  \begin{split}
    \left| \rho_k(\xi) - 1 \right| 
=\left| g(\xi-\xi_k, \sigma_k) - 1 \right| 
\leq e |\sigma_k| |\xi - \xi_k| \\
\leq e |\sigma_k| |\delta_{k-1}| 
\leq e \left( 2b k^{-5/2} \right) \left( k^{-3/2} \right) 
= 2 b e k^{-4},
  \end{split}
\end{equation}
for $\xi \in \left[\xi_{k+1}, \xi_{k-1}\right]$ and $k \geq k_0$.

Thus, for $\xi \in [\xi_{k+1}, \xi_k]$, 
using $0 \leq a_k(\xi) \leq 1$, (\ref{eq:rho-rhok}), (\ref{eq:akp}), (\ref{eq:akpp}), and (\ref{eq:rho-xi1})--(\ref{eq:rho-xi3}),
we obtain the following estimates:
\begin{equation}\label{eq:r-k1}
\begin{aligned}
\left| \rho(\xi) - 1 \right| 
&= \left|\rho_k(\xi) + (\rho_{k+1}(\xi) - \rho_{k}(\xi))a_k(\xi) - 1 \right| \\
&\leq \left|\rho_k(\xi) - 1 \right| + \left|\rho_{k+1}(\xi) - \rho_{k}(\xi) \right| \\
&\leq \left|\rho_k(\xi) - 1 \right| + \left|\rho_{k+1}(\xi) - 1\right| + \left|\rho_{k}(\xi) - 1 \right| \\
&\leq 6 b e k^{-4},
\end{aligned}
\end{equation}
\begin{equation}\label{eq:r-k2}
\begin{aligned}
\left| \rho'(\xi) \right| 
&= \left| \rho'_k(\xi) + (\rho'_{k+1}(\xi) - \rho'_k(\xi))a_k(\xi) + (\rho_{k+1}(\xi) - \rho_{k}(\xi))a_k'(\xi) \right| \\
&\leq \left| \rho'_k(\xi)\right| + \left| \rho'_{k+1}(\xi) - \rho'_k(\xi) \right| + \left| (\rho_{k+1}(\xi) - \rho_{k}(\xi))a_k'(\xi) \right| \\
&\leq 2ebk^{-5/2} + 4ebk^{-5/2} + (4ebk^{-4})(3c k^{3/2}) \\
&\leq (6eb + 12ceb)k^{-5/2},
\end{aligned}
\end{equation}
\begin{equation}\label{eq:r-k3}
\begin{aligned}
\left| \rho''(\xi) \right| 
&= \left| \rho''_k(\xi) + (\rho''_{k+1}(\xi) - \rho''_{k}(\xi))a_k(\xi) + 2(\rho'_{k+1}(\xi) - \rho'_{k}(\xi))a_k'(\xi) + (\rho_{k+1}(\xi) - \rho_{k}(\xi))a_k''(\xi) \right| \\
&\leq \left| \rho''_k(\xi)\right| + \left|\rho''_{k+1}(\xi) - \rho''_{k}(\xi) \right| + 2\left| (\rho'_{k+1}(\xi) - \rho'_{k}(\xi))a_k'(\xi) \right| + \left| (\rho_{k+1}(\xi) - \rho_{k}(\xi))a_k''(\xi) \right| \\
&\leq 2ebk^{-5/2} + 4ebk^{-5/2} + 2(4ebk^{-5/2})(3c k^{3/2}) + (4ebk^{-4})(5c k^{3}) \\
&\leq 2ebck^{-1} + 4ebck^{-1} + 24ebck^{-1} + 20ebck^{-1} 
= 50cebk^{-1}.
\end{aligned}
\end{equation}

From (\ref{eq:rhoatk1}) we have 
$$
\lim_{\xi \to 0^-} \rho(\xi) = 1, 
\quad 
\lim_{\xi \to 0^-} \rho'(\xi) =0, 
\quad 
\lim_{\xi \to 0^-} \rho''(\xi) =0,
$$
$$
\rho(0) = 1,
\quad 
\rho'(0^-) =0,
\quad
\rho''(0^-) = 0,
$$
where $\rho'(0^-)$ and $\rho''(0^-)$ denote the left-hand derivatives of $\rho(\xi)$ at $\xi = 0$.
On the other hand, from (\ref{eq:r-k1}) we have
$$
\lim_{\xi \to 0^+} \rho(\xi) = 1, 
\quad
\rho'(0^+) = \lim_{\xi \to 0^+} \frac{\rho(\xi)-1}{\xi} = 0;
$$
from (\ref{eq:r-k2}) we have
$$
\lim_{\xi \to 0^+} \rho'(\xi) = 0, 
\quad 
\rho''(0^+) = \lim_{\xi \to 0^+} \frac{\rho'(\xi)-0}{\xi} = 0;
$$
from (\ref{eq:r-k3}) we have 
$$
\lim_{\xi \to 0^+} \rho''(\xi) = 0.
$$
Therefore $\rho(\xi)$ is $C^2$ at $\xi=0$.

Moreover, we can compute that
$$
\lim_{\xi \to 0} \kappa(\xi) = \lim_{\xi \to 0} \frac{|\rho^2 + 2 \rho'^{2} - \rho \rho''|}{(\rho^2 + \rho'^{2})^{3/2}} = 1.
$$
Hence, there exists $k_1 > k_0$ such that the curvature $\kappa(\xi) > \frac{1}{2}$ for any $\xi \in(0, \xi_{k_1}]$. 
We remark that
on $(-\pi, 0]\cup [\xi_1,\pi]$, the curvature $\kappa(\xi)  \equiv 1$.

To complete our construction of the curve $\gamma$, 
we redefine $\rho_{k}(\xi)$ by setting $\rho_{k}(\xi) = 1$ for $k\leq k_1$. 
All previous estimates remain valid.
The resulting curve $\gamma$, given by the polar equation
\begin{equation}\label{eq:gamma-rho-ex}
r = \rho(\xi)
\end{equation}
is a $C^2$-smooth closed curve with everywhere positive curvature.

Lemma \ref{lem:Hal} is proved.

\begin{remark}
  It is an interesting observation that each time the curve $\gamma$ passes through $q_k$ for $k > k_1$, 
  it crosses from one side of the circle $S^1$ to the other, provided $\sigma_k \neq 0$. 
  As a consequence, there are infinitely many intersection points between $\gamma$ and the circle $S^1$ in a neighborhood of $\xi = 0$ in our construction.
\end{remark}

Using the construction developed above, we have proved Theorem \ref{thm:c2} in the case of $\mathbb{R}^3$. 
Next, we extend this construction to $\mathbb{R}^n$.

Let $K^2 \subset \mathbb{R}^3$ denote the $C^2$ convex cone constructed above, which admits a billiard trajectory (depending on a parameter $a$) with infinitely many reflections. When projected onto the plane $\mathcal{P}^2:=\{x \in \mathbb{R}^3 \mid x^3 = 1\}$, the trajectory has a unique limit point $q = (1, 0, 1)$ (independent of $a$).  
Near $q$, the cone $K^2$ is parameterized by
$$
r( x^2, t) = (t f_1(x^2), t x^2, t),
\quad t \in \mathbb{R}_{>0},
\quad x^2 \in (-1, 1), 
$$  
where $f_1$ is the $C^2$ function defining $\gamma$ in (\ref{eq:gamma-rho-ex}) for $\xi\in (-\frac{\pi}{2}, \frac{\pi}{2})$.
We can always take $k_1$ large enough in the construction of $\gamma$ such that $x^1 = f_1(x^2)$ coincides with the unit circle $(x^1)^2+(x^2)^2=1$ in $\mathcal{P}^2$ for $x^2 \in (-1, 0]\cup [\frac{1}{3},1)$.

Let  
$$
\begin{aligned}
  & e_1 := \frac{\partial r(x^2, t)}{\partial t} = (f_1(x^2), x^2, 1),\\
 & e_2 := \frac{1}{t} \frac{\partial r(x^2, t)}{\partial x^2} = (f_1'(x^2), 1, 0),
\end{aligned}
$$  
which form a basis of the tangent space of $T_{r(x^2,t)}K^2$.  
Let $l_k$, $k>k_1$ be the sequence of oriented lines associated to the trajectory constructed in the case of $\mathbb{R}^3$.
Let ${v}_{k}$, $\|{v}_{k}\| = 1$, be the direction of $l_k$.  
Since $l_k$ reflects to $l_{k+1}$ at $p_{k+1}$,  
we have, at the reflection point,  
\begin{equation}\label{eq:eofvinr3}
  \langle {{v}}_{k}, e_1 \rangle = \langle {{v}}_{k+1}, e_1 \rangle, 
\quad
\langle {{v}}_{k}, e_2 \rangle = \langle {{v}}_{k+1}, e_2 \rangle.
\end{equation}

Consider the hypersurface $\tilde{K}^{n-1}$ in $\mathbb{R}^n$ given by  
\begin{equation}\label{eq:rtilde}
  \begin{aligned}
    \tilde{r}_+(x^2, \ldots, x^{n-1}, t) = (t F_1(x^2, \ldots, x^{n-1}), t x^2, \ldots, t x^{n-1}, t),
\quad  
t>0,
\quad
(x^2, \ldots, x^{n-1}) \in D^{n-2}, \\
\tilde{r}_-(x^2, \ldots, x^{n-1}, t) = (-t \sqrt{1- \sum_{i=2}^{n-1}(x^i)^2}, t x^2, \ldots, t x^{n-1}, t),
\quad  
t>0,
\quad  
(x^2, \ldots, x^{n-1}) \in \overline {D}^{n-2},  
  \end{aligned}
\end{equation} 
where  
$$
F_1(x^2, \ldots, x^{n-1}) := \sqrt{f_1(x^2)^2 - (x^3)^2 - \cdots - (x^{n-1})^2},
$$  
$$
D^{n-2} := \{(x^2, \ldots, x^{n-1}) \in \mathbb{R}^{n-2} \mid (x^2)^2 + \cdots + (x^{n-1})^2 < 1\},$$ 
and $\overline {D}^{n-2}$ is the closure of $D^{n-2}$. Let $\tilde{\gamma}$ be the intersection of $\tilde{K}^{n-1}$ with $\mathcal{P}^{n-1}:=\{x \in \mathbb{R}^n \mid x^n = 1\}$. Then 
$\tilde{\gamma}$ is a $C^2$-smooth, strictly convex, closed submanifold of $\mathcal{P}^{n-1}$ with nondegenerate second fundamental form (see Lemma \ref{lem:hessian} below). 

Embed $\mathbb{R}^3$ in $\mathbb{R}^n$ as  
$$
(y^1, y^2, y^3) \mapsto (y^1, y^2, 0, \ldots, 0, y^3).
$$  
Let $\tilde{l}_k$, $\tilde{{v}}_{k}$, $\tilde{p}_k$, and $\tilde{q}$ be the images of $l_k$, ${{v}}_{k}$, $p_k$, and $q$  
under this embedding, respectively.  
Note that $\tilde{q}=(1,0, \ldots, 0, 1)$ lies in the part of $\tilde{K}^{n-1}$ parameterized by $\tilde{r}_+$, and $r(x^2, t)$ with $t>0$ and $x^2 \in (-1,1)$ is mapped to $\tilde{r}_+(x^2,0,\ldots, 0,t)$.

Let  
$$
\begin{aligned}
  & \tilde{e}_1 := \frac{\partial \tilde{r}_+(x^2, \ldots, x^{n-1}, t)}{\partial t} \bigg|_{x^3 = \ldots = x^{n-1} = 0},\\
  & \tilde{e}_j := \frac{1}{t} \frac{\partial \tilde{r}_+(x^2, \ldots, x^{n-1}, t)}{\partial x^j} \bigg|_{x^3 = \ldots = x^{n-1} = 0}, \quad 2 \leq j \leq n-1,
\end{aligned}
$$  
be the $(n-1)$ independent tangent vectors that form a basis of $T_{\tilde{r}_+(x^2,0,\ldots, 0,t)}\tilde{K}^{n-1}$. Then  
$$
\tilde{e}_1 = (f_1(x^2), x^2, 0, \ldots, 0, 1),
\quad 
\tilde{e}_2 = (f_1'(x^2), 1, 0, \ldots, 0, 0),
$$ 
$$
\tilde{e}_j = (0, \ldots, 0, 1 \text{ (in the $j$-th place)}, 0, \ldots, 0), \quad 3 \leq j \leq n-1.
$$  
Then, at the point $\tilde{p}_{k+1}$, we have  
$$
\langle \tilde{{v}}_{k}, \tilde{e}_j \rangle = 0, \quad \text{for } 3 \leq j \leq n-1, 
$$  
and from (\ref{eq:eofvinr3}) we have  
$$
\begin{aligned}
  \langle \tilde{{v}}_{k}, \tilde{e}_1 \rangle &=  
  \langle {{v}}_{k}, e_1 \rangle = \langle {{v}}_{k+1}, e_1 \rangle
  = \langle \tilde{{v}}_{k+1}, \tilde{e}_1 \rangle, \\
  \langle \tilde{{v}}_{k}, \tilde{e}_2 \rangle &= 
  \langle {{v}}_{k}, e_2 \rangle = \langle {{v}}_{k+1}, e_2 \rangle
  = \langle \tilde{{v}}_{k+1}, \tilde{e}_2 \rangle.
\end{aligned}
$$  
Thus, the lines $\tilde{l}_k$, $k > k_1$, form a billiard trajectory in $\tilde{K}^{n-1}$ with infinitely many reflections.

To complete the proof of Theorem \ref{thm:c2}, we prove the following lemma.
\begin{lemma}\label{lem:hessian}
  Let $\tilde{K}^{n-1}$ be defined by (\ref{eq:rtilde}). 
  Then 
  $\tilde{\gamma}:= \tilde{K}^{n-1}\cap \mathcal{P}^{n-1}$ is a $C^2$-smooth, strictly convex, closed submanifold of $\mathcal{P}^{n-1}$ with nondegenerate second fundamental form. 
\end{lemma}
\textbf{Proof.}
By the properties of $f_1(x^2)$, $\tilde{\gamma}$ is a closed $C^2$-smooth hypersurface in $\mathcal{P}^{n-1}$, and coincides with the unit sphere of $\mathcal{P}^{n-1}$ for $x^1 < f_1(\frac{1}{3})$. 

To prove that $\tilde{\gamma}$ is strictly convex with non-degenerate second fundamental form everywhere, it suffices to show that the Hessian matrix of $F_1$ is negative definite on $D^{n-2}$.

The Hessian matrix of $F_1(x^2, \ldots, x^{n-1})$ is given by
$$
\text{Hess}(F_1)
=
\left(\frac{\partial^2 F_1}{\partial x^i \partial x^j}\right)_{2\leq i,j\leq n-1},
$$
where for $i = 2$:
  $$
  \frac{\partial^2 F_1}{\partial x^2 \partial x^2} 
  = \frac{f_1(x^2)f_1''(x^2) + f_1'(x^2)^2}{F_1} - \frac{f_1(x^2)^2 f_1'(x^2)^2}{F_1^3},
  $$
for $i = 2$ and $j \geq 3$:
  $$
  \frac{\partial^2 F_1}{\partial x^2 \partial x^j} = \frac{f_1(x^2)f_1'(x^2) x^j}{F_1^3},
  $$
for $i \geq 3$:
  $$
  \frac{\partial^2 F_1}{\partial x^i \partial x^i} = -\frac{1}{F_1} - \frac{(x^i)^2}{F_1^3},
  $$
for $ i,j \geq 3 $ and $i \neq j$:
  $$
  \frac{\partial^2 F_1}{\partial x^i \partial x^j} = -\frac{x^i x^j}{F_1^3}.
  $$
Let $w= (w^2,\ldots, w^{n-1})^T \in \mathbb{R}^{n-2}$. 
We compute the quadratic form $w^T\text{Hess}(F_1)w$:
\begin{align*}
  \sum_{i,j=2}^{n-1}\frac{\partial^2 F_1}{\partial x^i \partial x^j} w^i w^j  
  &= \left( \frac{f_1(x^2)f_1''(x^2) + f_1'(x^2)^2}{F_1} 
  - \frac{f_1(x^2)^2 f_1'(x^2)^2}{F_1^3} \right) (w^2)^2 + \frac{2 f_1(x^2)f_1'(x^2)}{F_1^3} \sum_{j=3}^{n-1} x^j w^j w^2 \\
  &\quad - \sum_{i=3}^{n-1} \left( \frac{1}{F_1} + \frac{(x^i)^2}{F_1^3} \right) (w^i)^2 - \frac{2}{F_1^3} \sum_{3 \leq i < j\leq n-1} x^i x^j w^i w^j\\
  &= \frac{f_1(x^2)f_1''(x^2) + f_1'(x^2)^2}{F_1} (w^2)^2 
  - \frac{1}{F_1^3}\left( 
    f_1(x^2)^2 f_1'(x^2)^2 (w^2)^2 - 2 f_1(x^2) f_1'(x^2) w^2 \sum_{j=3}^{n-1} x^j w^j
   \right)\\
   &\quad 
  - \frac{\sum_{i=3}^{n-1} (w^i)^2}{F_1} - 
  \frac{\sum_{i=3}^{n-1} (x^i)^2(w^i)^2 +2\sum_{3 \leq i < j\leq n-1} x^i x^j w^i w^j }{F_1^2}\\
  &= \frac{f_1(x^2)f_1''(x^2) + f_1'(x^2)^2}{F_1} (w^2)^2 
  - \frac{1}{F_1^3}\left( 
    f_1(x^2) f_1'(x^2) w^2 -\sum_{j=3}^{n-1} x^j w^j
   \right)^2 - \frac{\sum_{i=3}^{n-1} (w^i)^2}{F_1}.
\end{align*}

To show $w^T\text{Hess}(F_1)w <0$ for any $w\neq 0$, it suffices to prove that $f_1(x^2)f_1''(x^2) + f_1'(x^2)^2<0$ for $x^2\in (-1,1)$.

For $x^2\in (-1,0]\cup [\frac{1}{3},1)$, where $f_1(x^2)=\sqrt{1-(x^2)^2}$, we have 
$$
f_1'(x^2) = \frac{-x^2}{\sqrt{1-(x^2)^2}},
\quad 
f_1''(x^2)= \frac{-1}{(1-(x^2)^2)^{3/2}}.
$$
Therefore 
$$
f_1(x^2)f_1''(x^2) + f_1'(x^2)^2= \sqrt{1-(x^2)^2}\frac{-1}{(1-(x^2)^2)^{3/2}} + \left(\frac{-x^2}{\sqrt{1-(x^2)^2}}\right)^2 = -1.
$$
For $x^2\in (0,\frac{1}{3})$, since the curvature $\kappa >\frac{1}{2}$ for the curve $\gamma$, we have 
$$
\kappa = \frac{|f_1''(x^2)|}{(1+f_1'(x^2)^2)^{3/2}} > \frac{1}{2}.
$$ 
Since $f_1''(x^2)<0$, we have 
$$
\frac{-f_1''(x^2)}{(1+f_1'(x^2)^2)^{3/2}} > \frac{1}{2},
$$ 
and therefore 
\begin{equation}\label{eq:f1pp}
  f_1''(x^2)< -\frac{1}{2}(1+f_1'(x^2)^2)^{3/2} < -\frac{1}{2}(1+f_1'(x^2)^2).
\end{equation}
Since $f_1'(0)=0$ and $f_1''(x^2)<0$, $f_1'(x^2)$ is decreasing on $(0,\frac{1}{3})$, hence
\begin{equation}\label{eq:f1p}
  0 > f_1'(x^2)> f_1'(\frac{1}{3}) = -\frac{1}{2\sqrt{2}}.
\end{equation}
Consequently, $f_1(x^2)$ is also decreasing:
\begin{equation}\label{eq:f1}
 1 > f_1(x^2) > f_1(\frac{1}{3})= \frac{2\sqrt{2}}{3}.
\end{equation}
Using (\ref{eq:f1pp})--(\ref{eq:f1}), we obtain
$$
f_1(x^2)f_1''(x^2) + f_1'(x^2)^2 
< f_1(x^2) \left(-\frac{1}{2}(1+f_1'(x^2)^2)\right) +  f_1'(x^2)^2 
< \frac{2\sqrt{2}}{3}\left(-\frac{1}{2}(1+f_1'(x^2)^2)\right)+  f_1'(x^2)^2 
$$
$$
= f_1'(x^2)^2 \left(1-\frac{\sqrt{2}}{3}\right) - \frac{\sqrt{2}}{3}
< \frac{1}{8}\left(1-\frac{\sqrt{2}}{3}\right)- \frac{\sqrt{2}}{3} <0.
$$
Therefore, the Hessian matrix of $F_1$ is negative definite on $D^{n-1}$, which shows that $\tilde{\gamma}$ is strictly convex with nondegenerate second fundamental form for $x^1>0$. Combined with the fact that $\tilde{\gamma}$ coincides with the unit sphere in $\mathcal{P}^{n-1}$ for $x^1< f_1(\frac{1}{3})$, we have proved the assertion in Lemma \ref{lem:hessian}.

\section{The Elliptic Cone in $\mathbb{R}^3$} 
  
\subsection{The first integral}
Let $K_e \subset \mathbb{R}^3$ be an elliptic cone defined by (\ref{eq:elliptic-cone}).
We have the following lemma.
\begin{lemma}\label{lem:int-ell}
	The Birkhoff billiard inside $K_e$ 
  admits the first intergral
	$$
  I_2 = 
  a^2 m_{2,3}^2 + b^2 m_{1,3}^2 - m_{1,2}^2.
	$$
\end{lemma}
\textbf{Proof.}
To show that $I_2$ is a first integral, we use the methods of \cite{DM}.
In \cite{DM}, differential equations on a billiard table admitting a first billiard integral were found.
Let a billiard table be parameterized by 
$ r = (r^1(u), r^2(u), r^3(u)) $, 
where $ u = (u^1, u^2)$.
The functions
$$
v = (v^1, v^2, v^3), \quad m = (m_{2,3}, m_{1,3}, m_{1,2}),
$$
and consequently any function of $v, m$
$$
H(v, m),
$$
are invariant along trajectories between reflections. On the other hand, the function
$$
G(s_1, s_2, u),
$$
where $G$ is some function and
$$
s_j := \left\langle v, \frac{\partial r}{\partial u^j} \right\rangle,
$$
is invariant at the moment of reflection. 
Hence, $H$ is a first billiard integral of the billiard table if there exists a function $G$ such that the identity
\begin{equation}\label{eq:method}
H(v, m) = G(s_1, s_2, u),
\quad 
\text{ where } x^i = r^i(u),
\end{equation}
holds for all $v$ with $ \|v\| = 1$. 
This identity (\ref{eq:method}) provides equations for $r^i(u)$ (see \cite{DM} and examples therein).

Let us parameterize the elliptic cone $K_e$ as 
$$
r(u^1, u^2) = (a u^1,b u^2, \sqrt{(u^1)^2 + (u^2)^2}).
$$
We have
$$
\frac{\partial r}{\partial u^1}(u^1, u^2) = (a, 0, \frac{u^1}{\sqrt{(u^1)^2 + (u^2)^2}}),
\quad 
\frac{\partial r}{\partial u^2}(u^1, u^2) = (0, b, \frac{u^2}{\sqrt{(u^1)^2 + (u^2)^2}}),
$$
$$
s_1 = a v^1 + \frac{u^1 v^3}{\sqrt{(u^1)^2 + (u^2)^2}}, 
\quad 
s_2 = b v^2 +  \frac{u^2 v^3}{\sqrt{(u^1)^2 + (u^2)^2}},
$$
$$
m_{1,2} = a u^1 v^2 - b u^2 v^1 , 
\quad 
m_{1,3} = a u^1 v^3 - \sqrt{(u^1)^2+ (u^2)^2} v^1 , 
\quad  
m_{2,3} = b u^2 v^3 - \sqrt{(u^1)^2+ (u^2)^2} v^2 .
$$
We can check that the equation 
$$
a^2 m_{2,3}^2 + b^2 m_{1,3}^2 -  m_{1,2}^2 = h_{1,1}(u) s_1^2 +  h_{2,2}(u) s_2^2  + h_{1,2}(u) s_1 s_2 + h_0(u) 
$$
holds when 
$$ 
h_{1,1} = -b^2 (u^1)^2 - (1+b^2) (u^2)^2 , 
\quad 
h_{2,2}= - (a^2+1) (u^1)^2 - a^2 (u^2)^2 ,
$$
$$
h_{1,2}= 2 u^1 u^2,  
\quad 
h_0(u)= b^2(a^2+1)(u^1)^2+ a^2(b^2+1) (u^2)^2 .
$$

Lemma \ref{lem:int-ell} is proved.

\begin{remark}
  It is remarkable that billiards in elliptic cones $K_e$ have two families of caustics.
  One is the family of spheres (see \cite{MY}), and the other is a family of elliptic cones. More precisely, let us fix the values of the integrals $I_1=c_1>0$, $I_2=c_2>0$. Then one can check that 
  the corresponding trajectories are tangent to the cone 
  $K_{\lambda}$ defined by
  $
  \frac{(x^3)^2}{1 -\lambda} = \frac{(x^1)^2}{a^2+\lambda} + \frac{(x^2)^2 }{b^2+\lambda},
  $ where  $\lambda = -\frac{c_2}{c_1}$.
\end{remark}

\subsection{Estimate of the number of reflections}
In this section, we establish Theorem \ref{thm:elliptic}, 
which provides an upper bound for the number of reflections in any billiard trajectory inside $K_e$. 
The proof will be based on the following lemma.

\begin{lemma}\label{lem:est-ell}
  Let $l$ be an oriented line intersecting $K_e$ at points $p_1$ and $p_2$, 
  satisfying (\ref{eq:i-c}).
  Then 
  $$
  \angle p_1 O p_2 
  > \arcsin 
  \frac{2 a b \sqrt{c_1 c_2}}{a^2(b^2+1)c_1 + (b^2+1)c_2}.
  $$
\end{lemma}
\textbf{Proof.}
To prove this result, we first express $\sin \angle p_1 O p_2$ in terms of $I_1$, $I_2$ and $m_{1,2}$, and then find its minimum value as $m_{1,2}$ varies.

Let 
$$ 
p_i = (a t_i \cos \xi_i, b t_i \sin \xi_i, t_i),
\quad  i=1,2,
\quad \xi_1\neq \xi_2 \in [0, 2\pi), 
\quad t_1, t_2 > 0.
$$
Let the direction of $l$ be given by 
$$
 v= \frac{p_2-p_1}{\| p_2-p_1 \|}.
$$
We have by straightforward computation 
$$
v= \frac{ \left(
  a t_2 \cos \xi_2-a t_1 \cos \xi_1  ,
  b t_2 \sin \xi_2-b t_1 \sin \xi_1 ,
  t_2 - t_1 \right)}
  {\| p_2-p_1 \|},
$$
\begin{equation}\label{eq:mijs}
m_{1,2} 
= \frac{a b t_1 t_2 \sin{\left(\xi_2 - \xi_1\right)}}
{\| p_2-p_1 \|},\quad
m_{1,3} 
= \frac{a t_1 t_2 \left(\cos{\xi_1} - \cos{\xi_2}\right)}
{\| p_2-p_1 \|},\quad 
m_{2,3} 
=\frac{b t_1 t_2 \left(\sin{\xi_1} - \sin{\xi_2}\right)}
{\| p_2-p_1 \|}.
\end{equation}
Then we have
\begin{equation*}
\begin{split}
\frac{m_{1,2}^2}{ a^2 m_{2,3}^2 + b^2 m_{1,3}^2} 
&= \frac{1 + \cos \left(\xi_2 - \xi_1\right)}{2}
\end{split}
\end{equation*}
\begin{equation*}
  \begin{split}
  \frac{a^2 m_{2,3}^2 - b^2 m_{1,3}^2}{ a^2 m_{2,3}^2 + b^2 m_{1,3}^2} 
  &= \cos \left(\xi_1 + \xi_2\right).
  \end{split}
  \end{equation*}
Thus we can write $\cos\left(\xi_2 - \xi_1\right) $ and $\cos\left(\xi_2 + \xi_1\right) $
in terms of $I_1, I_2,$ and $m_{1,2}$
\begin{align}
  \label{eq:angle-to-integral1}
  \cos\left(\xi_2 - \xi_1\right) 
  =&  \frac{  2 m_{1,2}^2}{ a^2 m_{2,3}^2 + b^2 m_{1,3}^2} -1 = \frac{2 m_{1,2}^2}{m_{1,2}^2 + I_2}-1,\\
  \label{eq:angle-to-integral2}
  \cos\left(\xi_2 + \xi_1\right) 
  =& \frac{  a^2 m_{2,3}^2 - b^2 m_{1,3}^2}{ a^2 m_{2,3}^2 + b^2 m_{1,3}^2}  
=\frac{(a^2 + b^2)I_2 - 2 a^2 b^2 I_1 + (a^2+b^2 + 2a^2b^2) m_{1,2}^2}{(a^2- b^2 )(I_2 + m_{1,2}^2)}.
\end{align}

From Lemma \ref{lem:int-dis} we know 
$$
\text{dist}(l, O) = \sqrt{I_1}.
$$
Thus from (\ref{eq:area-p}) we obtain
\begin{equation}\label{eq:e1}
  \sin \angle p_1 O p_2 
 = \frac{\sqrt{I_1} \left\|p_2-p_1\right\| }{ \left\|p_1\right\| \left\| p_2\right\|}.
\end{equation}
To compute $\left\|p_2-p_1\right\|$,
we have from (\ref{eq:mijs}) that 
$$
\|p_2-p_1\|^2 = \frac{a^2 b^2 t_1^2 t_2^2 \sin^2\left(\xi_2 - \xi_1\right)}{m_{1,2}^2} = 
\frac{a^2 b^2 t_1^2 t_2^2 \left(1 - \cos^2\left(\xi_2 - \xi_1\right)\right)}{m_{1,2}^2}.
$$
Substituting (\ref{eq:angle-to-integral1}), we obtain
\begin{equation}\label{eq:e2}
\|p_2-p_1\|^2 = \frac{a^2 b^2 t_1^2 t_2^2}{m_{1,2}^2} \left( 1 - \left( \frac{2 m_{1,2}^2}{m_{1,2}^2 + I_2} - 1 \right)^2 \right)
= \frac{4 a^2 b^2 t_1^2 t_2^2 I_2}{\left(m_{1,2}^2 + I_2\right)^2}.
\end{equation}
To compute $\left\|p_1\right\|\left\|p_2\right\|$, we have 
\begin{equation*}
  \begin{split}
    \left\|p_1\right\|^2 \left\| p_2\right\|^2 
= & \ 
t_1^2 t_2^2 \left(1+a^2\cos^2 \xi_1 + b^2 \sin^2 \xi_1\right) \left(1+a^2 \cos^2 \xi_2 + b^2 \sin^2 \xi_2\right) 
\\
= & \ 
t_1^2 t_2^2 \left(1+ \frac{a^2 + b^2 }{2} + \frac{a^2 - b^2 }{2}\cos 2\xi_1 \right)\left(1+ \frac{a^2 + b^2 }{2} + \frac{a^2 - b^2 }{2}\cos 2\xi_2 \right)\\
= & \  
t_1^2 t_2^2\left((1+a^2)(1+b^2) + \frac{(2+ a^2 + b^2)(a^2 - b^2) }{2}\cos\left(\xi_2 - \xi_1\right)	\cos\left(\xi_2 + \xi_1\right) \right.
\\
& \left. + \frac{(a^2 - b^2)^2 }{4} \left(\cos\left(\xi_2 - \xi_1\right)^2 + \cos\left(\xi_2 + \xi_1\right)^2\right)\right).
  \end{split}
\end{equation*}
Substituting (\ref{eq:angle-to-integral1}), (\ref{eq:angle-to-integral2}), we get
\begin{equation}\label{eq:e3}
\left\|p_1\right\|^2 \left\| p_2\right\|^2 = t_1^2 t_2^2
\frac
{4 a^2 b^2 I_1 I_2 + \left( (1+a^2)(1+b^2) m_{1,2}^2 - ( a^2 b^2 I_1-I_2) \right)^2}
{\left(m_{1,2}^2 + I_2\right)^2}.
\end{equation}
Combining (\ref{eq:e1}), (\ref{eq:e2}), (\ref{eq:e3}), we obtain 
$$
\sin^2\angle p_1 O p_2  =  \frac {4 a^2 b^2 I_1 I_2}
{4 a^2 b^2 I_1 I_2 + \left( (1+a^2)(1+b^2) m_{1,2}^2 - ( a^2 b^2 I_1-I_2) \right)^2}.
$$

Since $a^2>b^2$, we have
$$
a^2 I_1-I_2 = (a^2-b^2) m_{1,3}^2 + (a^2 +1 )m_{1,2}^2 \geq (a^2 +1 )m_{1,2}^2,
$$
which implies
$$
m_{1,2}^2 \leq \frac{a^2 I_1-I_2}{a^2 +1}.
$$
Therefore, when $I_1=c_1>0$ and $I_2=c_2>0$ are fixed, $m_{1,2}^2$ is bounded by
$$
0 \leq m_{1,2}^2 \leq  \frac{a^2c_1 -c_2}{a^2 +1}.
$$

Define a function
$$
h(s) :=  \frac {4 a^2 b^2 c_1 c_2}{4 a^2 b^2 c_1 c_2 + \left( (1+a^2)(1+b^2)s - ( a^2 b^2 c_1-c_2) \right)^2}, \quad
s\in J:=\left[0, \frac{a^2c_1 -c_2}{a^2 +1}\right].
$$
Then
\begin{equation*}
\min_{s\in J} h(s) 
=
\frac {4 a^2 b^2 c_1 c_2}{4 a^2 b^2 c_1 c_2 + \max_{s\in J} f(s)^2},
\end{equation*}
where $f(s):=(1+a^2)(1+b^2)s - (a^2 b^2 c_1-c_2)$. Since $f(s)$ is linear in $s$, 
$\max_{s\in J} f(s)^2$ can only be achieved at the endpoints of $J$, i.e.,
$$
\max_{s\in J} f(s)^2 = \max\{f(0)^2, f(\frac{a^2c_1 -c_2}{a^2 +1})^2\}.
$$
Therefore,
\begin{equation*}
\begin{split}
\min_{s\in J} h(s) 
&= \min \{ h(0), h(\frac{a^2c_1-c_2}{a^2 +1})\} \\
&= \min \{ \frac{4 a^2 b^2 c_1 c_2}{(a^2 b^2 c_1 + c_2)^2}, \frac{4 a^2 b^2 c_1 c_2}{(a^2 c_1 + b^2 c_2)^2}\}.
\end{split}
\end{equation*}
Hence
$$
\sin \angle p_1 O p_2 
 \geq \min
 \{ 
  \frac{2 a b \sqrt{c_1 c_2}}{a^2 b^2 c_1 + c_2}, 
 \frac{2 a b \sqrt{c_1 c_2}}{a^2 c_1 + b^2 c_2} 
 \}
>
\frac{2 a b \sqrt{c_1 c_2}}{a^2(b^2+1)c_1 + (b^2+1)c_2}.
$$
Therefore
$$
  \angle p_1 O p_2 
  > \arcsin 
  \frac{2 a b \sqrt{c_1 c_2}}{a^2(b^2+1)c_1 + (b^2+1)c_2}.
$$

Lemma \ref{lem:est-ell} is proved.

\vspace*{1em}

\paragraph{Proof of Theorem \ref{thm:elliptic}.}
Consider a billiard trajectory inside $K_e$ consisting of a sequence (finite or infinite) of oriented lines $l_1, l_2, \ldots$;
for $k \geq 1$, $l_k$ passes through $p_k$ and $p_{k+1}$.
We count the number of reflections starting from $p_1$ (note that by convexity of $K_e$, there exists an oriented line $l_0$ that reflects at $p_1$ into $l_1$). 
When the trajectory reaches the point $p_N$, the number of reflections (including $p_1, \ldots, p_N$) is counted as $N$.

By (\ref{eq:alpha-theta}) we have
$$
\alpha_{k+1} = \alpha_{k} - \theta_k, 
\quad  k \geq 1.
$$ 
Thus, if $p_N$ lies on the trajectory, then
$$
\sum_{k=1}^{N-1}  \theta_k = \alpha_{1} - \alpha_{N} < \pi.
$$
From Lemma \ref{lem:est-ell} we have
$$
\theta_k 
> 
\arcsin 
\frac{2ab\sqrt{c_1 c_2}}
{a^2(b^2 + 1)c_1 + (b^2 + 1)c_2},
\quad 1\leq k \leq N-1.
$$
Therefore
$$
\pi > \sum_{k=1}^{N-1} \theta_k 
> 
(N-1) 
\arcsin 
\frac{2ab\sqrt{c_1 c_2}}
{a^2(b^2 + 1)c_1 + (b^2 + 1)c_2},
$$
which implies 
$$
N < \frac{\pi}{\arcsin \frac{2ab\sqrt{c_1 c_2}}{a^2(b^2 + 1)c_1 + (b^2 + 1)c_2}} + 1.
$$
Thus, the trajectory has at most 
$$
N_{c_1, c_2} = \left\lceil \frac{\pi }{\arcsin \frac{2ab\sqrt{c_1 c_2}}{ a^2(b^2 + 1)c_1 + (b^2 + 1)c_2}} \right\rceil 
$$
reflections.
It follows immediately that the number of reflections in any trajectory inside $K_e$ is finite. 
This provides an alternative proof of Theorem 2 in \cite{MY} for elliptic cones in $\mathbb{R}^3$.

This completes the proof of Theorem \ref{thm:elliptic}.

\vspace{2cm}

\noindent
{Andrey E. Mironov}\\
Sobolev Institute of Mathematics, Novosibirsk, Russia\\
Email: \texttt{mironov@math.nsc.ru}

\medskip
\noindent
{Siyao Yin}\\
Sobolev Institute of Mathematics, Novosibirsk, Russia\\
Email: \texttt{yinsiyao@outlook.com}

\end{document}